\documentclass[11pt]{extarticle}

\usepackage{amsmath}
\usepackage{amssymb}
\usepackage{amsthm}
\numberwithin{equation}{section}
\usepackage{authblk}
\usepackage{caption}
\usepackage{comment}
\usepackage{dsfont}
\usepackage[dvipsnames]{xcolor}
\usepackage{eurosym}
\usepackage{geometry}
\usepackage{lipsum}
\usepackage{mathrsfs}
\usepackage{mathtools}

\usepackage{natbib}
\usepackage[colorlinks=true, allcolors=blue]{hyperref}
\usepackage[OT1]{fontenc}
\usepackage{soul}
\usepackage{upquote}

\newtheorem{assumption}{Assumption}[section]

\newtheorem{lemma}{Lemma}[section]
\newtheorem{proposition}{Proposition}[section]
\newtheorem{remark}{Remark}[section]
\newtheorem{theorem}{Theorem}[section]
\newtheorem{corollary}{Corollary}[section]
\theoremstyle{remark}

\newtheorem{example}{Example}[section]

\DeclareMathOperator{\Erw}{\mathds E}

\DeclareMathOperator{\Prob}{\mathds P}

\begin{document}

\title{Estimating invertible processes in Hilbert spaces, with applications to functional ARMA processes}

\author[1]{Sebastian Kühnert}
\author[2]{Gregory Rice}
\author[3]{Alexander Aue}

\affil[1]{Department of Mathematics, Ruhr University Bochum, 44780 Bochum, Germany, \textsuperscript{*}\href{mailto:sebastian.kuehnert@ruhr-uni-bochum.de}{sebastian.kuehnert@ruhr-uni-bochum.de}}
\affil[2]{Department of Statistics and Actuarial Science, University of Waterloo, ON N2L 0A4 Waterloo, Canada, \textsuperscript{\dag}\href{mailto:grice@uwaterloo.ca}{grice@uwaterloo.ca}}
\affil[3]{Department of Statistics, University of California, Davis, CA 95616 Davis, USA, \textsuperscript{\ddag}\href{mailto:aaue@ucdavis.edu}{aaue@ucdavis.edu}}

\date{Jul 25, 2025}

\maketitle

\begin{abstract} Invertible processes are central to functional time series analysis, making the estimation of their defining operators a key problem. While asymptotic error bounds have been established for specific ARMA models on $L^2[0,1]$, a general theoretical framework has not yet been considered. This paper fills in this gap by deriving consistent estimators for the operators characterizing the invertible representation of a functional time series with white noise innovations in a general separable Hilbert space. Under mild conditions covering a broad class of functional time series, we establish explicit asymptotic error bounds, with rates determined by operator smoothness and eigenvalue decay. These results further provide consistency-rate estimates for operators in Hilbert space-valued causal linear processes, including functional MA, AR, and ARMA models of arbitrary order.
\end{abstract}

\noindent{\small \textit{MSC 2020 subject classifications:} 60G10, 62F12, 62R10}

\noindent{\small \textit{Keywords:} ARMA; functional time series; invertible processes; linear processes}

\section{Introduction}\label{Section: Introduction}
Over the past three decades, the field of \emph{functional data analysis} (fDA) has grown significantly, driven by increasing interest in high-dimensional data obtained from continuous observations across time, space, and frequency; see \cite{HorvathKokoszka2012} and \cite{RamsaySilverman2005} for textbook length reviews. In many cases, functional data are collected over time, such as continuous electricity price observations transformed into daily curves \citep{Liebl2013}. Such data are known as \emph{functional time series} (fTS). Recent advances in their analysis and modeling are discussed in \cite{Bosq2000} and \cite{kokoszka:2017:FDA-book}. A typical setup in fTS analysis (fTSA) assumes data reside in a separable, infinite-dimensional Hilbert space $\mathcal{H}$, which we consider throughout. Similar to scalar time series, linear models are most commonly used in fTSA, including $\mathcal{H}$-valued autoregressive (fAR) processes, functional moving average (fMA), and fARMA processes. Under regularity conditions, such processes admit stationary solutions, often represented as $\mathcal{H}$-valued linear and invertible processes. For testing stationarity in general fTS, see \cite{HorvathKokoszkaRice2014, AueVanDelft2020}, and \cite{vanDelftCharaciejusDette2021}. Invertible representations of stationary processes also arise in non-linear fTS models, such as $\mathcal{H}$-valued volatility processes within f(G)ARCH models, see \cite{Hoermannetal2013, Aueetal2017, Ceroveckietal2019, Kuehnert2020}.  
Early works such as \cite{Bosq2000} proposed consistent estimators for fAR operators. Limit theorems for $\mathcal{H}$-valued linear processes are available in \cite{MERLEVEDE1996, Merlevedeetal1997, Dueker2018}, and \cite{RackauskasSuquet2010}. \cite{AueKlepsch2017} and \cite{klepsch2017innovations} focused on point-wise estimation of operators in the linear and inverted representations, and \cite{Kuehnert2020} estimated operators of invertible linear processes originating from fGARCH models in $L^2[0,1]$. For fAR operators, \cite{CaponeraPanaretos2022} derived explicit rates under mild conditions, while \cite{KaraTerkiMourid2024} established exponential bounds and convergence rates using sieve estimates. 

The fMA and fARMA cases have received less attention. \cite{TurbillonETAL2008} estimated the fMA$(1)$ operator assuming commutativity, and \cite{kuenzer:2024} derived consistent estimators for fARMA operators using a two-step approach, involving sequential AR and MA estimation with dimension reduction via principal components. Central limit theorems for fAR operator estimates are derived in \cite{mas2009linear}.

This article formulates consistency results with quantifiable rates for infinite-dimensional operators defining general $\mathcal{H}$-valued invertible processes with white noise innovations, reyling on Yule--Walker equations and a Tikhonov-regularized estimator, avoiding assumptions tied to specific models. The consistency results are established under a Sobolev-type condition and assuming general consistency rates for certain lagged covariance operator estimates that hold in many standard situations, such as under $L^p$-$m$-approximability, a weak dependence concept for fTSA  introduced in \cite{hormann:kokoszka:2010}. We apply these findings to derive consistency for the infinite-dimensional operators defining causal linear processes, including fAR, fMA, and fARMA processes with arbitrary AR and MA orders. Notably, our results are derived under mild conditions, permitting multi-dimensional eigenspaces of the covariance operator, and without requiring specific rates of eigenvalue decay. Provided the operators exhibit a sufficiently high degree of smoothness and the eigenvalue decay of the error covariance operators is slow, our approach nearly achieves the parametric rate $N^{-1/2}.$

The paper is organized as follows: Section \ref{Section: Preliminaries} introduces notation and fundamental concepts. Section \ref{Section: Estimation in the invertible and linear representation} derives estimators for operators of invertible processes. Sections \ref{Section: Linear process representation}--\ref{Section: Implications to ARMA} apply these results to obtain consistent estimators for operators in linear and fARMA processes. Section \ref{Section: Summary} concludes, with proofs in Appendix \ref{Section Proofs}. 

\section{Preliminaries}\label{Section: Preliminaries}

\subsection{Notation}

The additive identity in linear spaces is denoted by $0$, and the identity map by $\mathbb{I}$. In Cartesian product spaces $V^n$, $n \in \mathbb{N}$, scalar multiplication and addition are defined component-wise. Let $(\mathcal{B}, \|\cdot\|)$ be a Banach space and $(\mathcal{H}, \langle\cdot, \cdot\rangle)$, $(\mathcal{H}_\star, \langle\cdot, \cdot\rangle_\star)$ Hilbert spaces with their respective norms and inner products. Then, $\mathcal{B}^n$ with norm $\|x\|^2 = \sum_{i=1}^n\|x_i\|^2$ is a Banach space, and $\mathcal{H}^n$ with inner product $\langle x, y\rangle = \sum_{i=1}^n\langle x_i, y_i\rangle$ is a Hilbert space. The spaces of bounded linear, Hilbert--Schmidt (H-S), and nuclear (trace-class) operators $A\colon \mathcal{H} \to \mathcal{H}_\star$ are denoted by $\mathcal{L}_{\mathcal{H}, \mathcal{H}_\star}$, $\mathcal{S}_{\mathcal{H}, \mathcal{H}_\star}$, and $\mathcal{N}_{\mathcal{H}, \mathcal{H}_\star}$, respectively, with norms $\|\cdot\|_{\mathcal{L}}$, $\|\cdot\|_{\mathcal{S}}$, and $\|\cdot\|_{\mathcal{N}}$, and H-S inner product $\langle\cdot,\cdot \rangle_{\mathcal{S}}$. For $\mathcal{T} \in \{\mathcal{L,S,N}\}$, set $\mathcal{T_H} \coloneqq \mathcal{T}_{\mathcal{H,H}}$. The projection onto the closed subspace spanned by $a_m, \dots, a_n$ of a \emph{complete orthonormal system} (CONS) $(a_j)_{j\in\mathbb{N}}$ in $\mathcal{H}$ is denoted $\coprod^{a_n}_{a_m}\colon \mathcal{H} \to \mathrm{span}\{a_m, \dots, a_n\} \subset \mathcal{H}$. All random variables are defined on a common probability space $(\Omega, \mathfrak{A}, \mathbb{P})$. For a sequence $(c_n)_{n\in\mathbb{N}} \subset (0,\infty)$, $X_n = \mathrm{o}_{\mathbb{P}}(c_n)$ and $X_n = \mathrm{O}_{\mathbb{P}}(c_n)$ as $n \to \infty$ denote $X_n / c_n$ converges in probability and is stochastically bounded, respectively. The space of random variables $X \in \mathcal{H}$ with $\mathbb{E}\|X\|^p < \infty,$ $p \in [1,\infty)$, is denoted by $L^p_{\mathcal{H}} = L^p_{\mathcal{H}}(\Omega, \mathfrak{A}, \mathbb{P})$. A sequence $(X_k) \subset L^1_\mathcal{H}$ is \emph{centered} if $\mathbb{E}(X_k) = 0$ for all $k$ with expectation understood as a Bochner integral. For sequences $(a_n)_{n\in\mathbb{N}}, (b_n)_{n\in\mathbb{N}} \subset \mathbb{R}$, $a_n \asymp b_n$ as $n \to \infty$ indicates asymptotic equivalence up to constant factors. Finally, $\mathds{1}_A(x)$ denotes the indicator function of a set $A,$ and $\delta_{ij}$ the \emph{Kronecker delta}. Unless stated otherwise, asymptotic statements refer to $N \to \infty$.  

\subsection{Weak stationarity and lagged  (cross-)covariance operators}

Let $(\mathcal{H}, \langle\cdot,\cdot\rangle)$ and $(\mathcal{H}_\star, \langle\cdot,\cdot\rangle_\star)$ be separable Hilbert spaces with norms $\|\cdot\| = \langle\cdot,\cdot\rangle^{1/2}$ and $\|\cdot\|_\star = \langle\cdot,\cdot\rangle_\star^{1/2}$, respectively. The \emph{cross-covariance operator} of two random variables $X \in L^2_{\mathcal{H}}$ and $Y \in L^2_{\mathcal{H}_\star}$ is defined by $\mathscr{C}_{\!X,Y} \coloneqq \mathbb{E}[X-\mathbb{E}(X)] \otimes [Y-\mathbb{E}(Y)],$ where $x \otimes y \coloneqq \langle x, \cdot\rangle y$ denotes the \emph{tensor product}. The operator $\mathscr{C}_{\!X,Y} \in \mathcal{N}_{\mathcal{H}, \mathcal{H}_\star}$ has the adjoint $\mathscr{C}^\ast_{\!X,Y} = \mathscr{C}_{Y,X} \in \mathcal{N}_{\mathcal{H}_\star, \mathcal{H}}$, and 
\begin{align}\label{Nuclear norm cross-cov op}
  \|\mathscr{C}_{\!X,Y}\|_{\mathcal{N}} \leq \Erw\!\|X-\Erw(X)\|\|Y-\Erw(Y)\|_\star\,.    
\end{align}
The \emph{covariance operator} of $X \in L^2_{\mathcal{H}}$ is $\mathscr{C}_{\!X} \coloneqq \mathscr{C}_{\!X,X}.$ It is positive semi-definite, self-adjoint, and 
\begin{align}\label{Nuclear norm cov op}
  \|\mathscr{C}_{\!X}\|_{\mathcal{N}} = \Erw\!\|X-\Erw(X)\|^2.    
\end{align}

A process $(X_k)_{k\in\mathbb{Z}} \subset L^2_{\mathcal{H}}$ is \emph{weakly stationary} if there exists $\mu \in \mathcal{H}$ such that $\mathbb{E}(X_k) = \mu$ for all $k$, and $\mathscr{C}_{X_k, X_\ell} = \mathscr{C}_{X_{k+h}, X_{\ell+h}}$ for all $h, k, \ell$. Such a process is a \emph{white noise} (WN) if $\mu = 0$, $\mathbb{E}\|X_k\|^2 > 0$ for all $k$, and $\mathscr{C}_{X_k, X_\ell} = 0$ for $k \neq \ell$. The \emph{lag-$h$ covariance operators} of a weakly stationary process $\boldsymbol{X} = (X_k)_{k \in \mathbb{Z}}$ are $\mathscr{C}^h_{\!\boldsymbol{X}} \coloneqq \mathscr{C}_{\!X_0, X_h}$ for $h \in \mathbb{Z}$, and the covariance operator is $\mathscr{C}_{\!\boldsymbol{X}} = \mathscr{C}^0_{\!\boldsymbol{X}}$. For two processes $\boldsymbol{X} = (X_k)_{k \in \mathbb{Z}} \subset L^2_{\mathcal{H}}$ and $\boldsymbol{Y} = (Y_k)_{k \in \mathbb{Z}} \subset L^2_{\mathcal{H}_\star}$, the \emph{lag-$h$ cross-covariance operators} are defined by $\mathscr{C}^h_{\!\boldsymbol{X}, \boldsymbol{Y}} \coloneqq \mathscr{C}_{\!X_0, Y_h}$ for $h \in \mathbb{Z}$, provided $\mathscr{C}_{X_k, Y_\ell} = \mathscr{C}_{X_{k+h}, Y_{\ell+h}}$ for all $h, k, \ell$.

\section{Estimation of the operators in  invertible processes}\label{Section: Estimation in the invertible and linear representation}

This section develops an estimation procedure for operators of \emph{invertible} processes $(X_k)_{k\in\mathbb{Z}}$ in a separable Hilbert space $(\mathcal{H}, \langle\cdot, \cdot\rangle).$ While such processes are well-defined for bounded linear operators, our approach requires H-S  operators. Throughout, we assume the following. 

\begin{assumption}\label{As: Features of (eps_k)} 
$\boldsymbol{\varepsilon}=(\varepsilon_k)_{k\in\mathbb{Z}}\subset \mathcal{H}$ is a WN with injective covariance operator $\mathscr{C}_{\boldsymbol{\varepsilon}}\colon\mathcal{H}\rightarrow\mathcal{H}.$
\end{assumption}

We call a centered process $(X_k)_{k\in\mathbb{Z}}\subset\mathcal{H}$ \emph{invertible} if it is weakly stationary and satisfies
\begin{align}\label{identity of invertible process}
  X_k =  \sum^\infty_{j=1}\psi_j(X_{k-j}) + \varepsilon_k,\quad k\in\mathbb{Z},
\end{align}
where $(\varepsilon_k)_{k\in\mathbb{Z}}$ is a WN, and $\psi_j\in\mathcal{L_H}$ are bounded  linear operators with $\sum^\infty_{j=1}\|\psi_j\|_\mathcal{L} < \infty,$ ensuring the series converges a.s. 
 \frenchspacing and in $L^2_\mathcal{H}.$   Such processes arise in various contexts, including invertible linear processes \citep[Chapter 7]{Bosq2000}, encompassing fAR, as well as fMA and fARMA models under natural conditions \citep{Spangenberg2013, kuenzer:2024}. Moreover, estimating fAR(MA) parameters in invertible processes is relevant for f(G)ARCH models, as their point-wise squared counterparts correspond to fAR(MA) processes \citep{Hoermannetal2013, Aueetal2017, Kuehnert2020} with weak WN innovations.

\begin{assumption}\label{As: Inv process} $(X_k)_{k\in\mathbb{Z}}\subset L^2_\mathcal{H}$ is a centered, invertible process with representation \eqref{identity of  invertible process}, with operators $\psi_j \in  \mathcal{S}_{\mathcal{H}}$.
\end{assumption}

Consistent estimates for the operators $\psi_j$ are derived based on a sample $X_1, \dots, X_N$ from ${\boldsymbol{X}} = (X_k)$. The proposed estimators are based on lagged (cross-)covariance operators of processes
\begin{align}\label{process X in Cartesian space}
  X^{[L]}_k \coloneqq (X_k, X_{k-1}, \dots, X_{k-L+1})^\top, \quad k \in \mathbb{Z},    
\end{align}
in the Cartesian product Hilbert space $\mathcal{H}^L$. The lag parameter $L = L_N \in \mathbb{N}$ is assumed to grow with the sample size. We define the empirical covariance operator $\hat{\mathscr{C}}_{\!\boldsymbol{X}^{[L]}}$ of $\boldsymbol{X}^{[L]} = (X^{[L]}_k)_k$ and the empirical lag-1 cross-covariance operator $\hat{\mathscr{D}}_{\!\boldsymbol{X}^{[L]}, \boldsymbol{X}}$ of $\boldsymbol{X}^{[L]}$ and $\boldsymbol{X}$ by
\begin{alignat}{2}\label{empirical lag-h-covariance operators}
   \hat{\mathscr{C}}_{\!\boldsymbol{X}^{[L]}} &\coloneqq \frac{1}{N} \sum_{k=L}^N X^{[L]}_k \otimes X^{[L]}_k, \quad &L \leq N,\\  \hat{\mathscr{D}}_{\!\boldsymbol{X}^{[L]}, \boldsymbol{X}} &\coloneqq \frac{1}{N} \sum_{k=L}^{N-1} X^{[L]}_k \otimes X_{k+1}, \quad &L < N. \label{empirical lag-1-cross covariance operator}
\end{alignat}
To establish consistency for the estimates for $\psi_j$, we assume the following.
\begin{assumption}\label{As: cons rate for lagged cov and cross-cov ops} For some  sequence $L=L_N\rightarrow\infty$ with $L=\mathrm o(N^{1/3})$ holds
\begin{align}
  \|\hat{\mathscr{C}}_{\!\boldsymbol{X}^{[L]}} - \mathscr{C}_{\!\boldsymbol{X}^{[L]}}\|_{\mathcal{S}} &= \mathrm O_{\Prob}(L^{3/2}N^{-1/2}), \label{up b cov}\\
\|\hat{\mathscr{D}}_{\!\boldsymbol{X}^{[L]}, \boldsymbol{X}} - \mathscr{D}_{\!\boldsymbol{X}^{[L]}, \boldsymbol{X}}\|_{\mathcal{S}} &= \mathrm O_{\Prob}(LN^{-1/2}). \label{up b cross-cov}
	\end{align}
\end{assumption}

\begin{remark} Assumption \ref{As: cons rate for lagged cov and cross-cov ops} holds for centered processes under $L^4$-$m$-approximability \citep{hormann:kokoszka:2010}, as shown in Corollary 3.1 of \cite{Kuehnert2025}, and for \emph{$\pi$-dependence}, a notion of weak dependence in metric spaces put forth by \cite{Kutta2025}. It also applies to i.i.d. processes with finite fourth moments \citep{HorvathKokoszka2012}. Notably, there are non-causal processes that satisfy Assumption \ref{As: cons rate for lagged cov and cross-cov ops}, such as \emph{linear processes} in $\mathcal{H}$, that is
\begin{align}\label{Definition LP}
X_k = \sum_{\ell=-\infty}^\infty \phi_\ell(\varepsilon_{k-\ell}),
\end{align}
where $\phi_\ell \in \mathcal{L_H}$ with $\sum^\infty_{\ell=-\infty} \|\phi_\ell\|_{\mathcal{L}} < \infty.$ Then,  Assumption \ref{As: cons rate for lagged cov and cross-cov ops} is satisfied if the linear process operators fulfill the stronger condition
$\sum^\infty_{k=-\infty}|k|\,\|\phi_k\|_{\mathcal{L}} < \infty.$
\end{remark}
\noindent The estimates for the operators $\psi_j$ in the inverted representation \eqref{identity of invertible process} are derived based on the approximate Yule--Walker equation
\begin{align}\label{approximate YW equation}
  \mathscr{D}_{\!\boldsymbol{X}^{[L]}, \boldsymbol{X}} = \Psi_{\!L}\mathscr{C}_{\!\boldsymbol{X}^{[L]}} + \sum_{\ell>L}\,\psi_\ell\mathscr{C}^{1-\ell}_{\!\!\boldsymbol{X}^{[L]}, \boldsymbol{X}},
\end{align}
\noindent where the process $\boldsymbol{X}^{[L]} = (X^{[L]}_k)_k\subset\mathcal{H}^L$ is defined in \eqref{process X in Cartesian space}, the operator $\Psi_{\!L} \in \mathcal{L}_{\mathcal{H}^L\!,\mathcal{H}}$ is given by
\begin{align}\label{Definition of operator-valued 1xL-matrix}
\Psi_{\!L} \coloneqq (\psi_1\,\cdots\,\psi_L),
\end{align}
and the series in \eqref{approximate YW equation} exists for all fixed $L\in\mathbb{N},$ as it converges in the H-S norm topology under Assumption \ref{As: Inv process}, see \eqref{inequ sum operators}. Note that $\Psi_{\!L}$ is identifiable from $\Psi_{\!L}\mathscr{C}_{\!\boldsymbol{X}^{[L]}}$ if $\mathscr{C}_{\!\boldsymbol{X}^{[L]}}$ is injective.
\begin{proposition}\label{Prop: Injectivity of C^L_X} 
Let Assumptions \ref{As: Features of (eps_k)}--\ref{As: Inv process} hold. Then, $\mathscr{C}_{\!\boldsymbol{X}^{[L]}}$ is injective for all $L\in\mathbb{N}.$
\end{proposition}
Given a sample $X_1, \dots, X_N$ from $\boldsymbol{X}$, we estimate $\Psi_{\!L}$ via Tikhonov regularization, which adds a damping term to the eigenvalues of $\mathscr{C}_{\!\boldsymbol{X}^{[L]}}$ \citep[see, e.g.,][]{hall:horo:2005}. This yields the estimator 
\begin{align}\label{Operator estimate for big Psi_L}
\hat{\Psi}_{\!L}\coloneqq\hat{\mathscr{D}}_{\!\boldsymbol{X}^{[L]}, \boldsymbol{X}}\hat{\mathscr{C}}^\dagger_{\!\!\boldsymbol{X}^{[L]}}\!\coprod^{\hat{c}_K}_{\hat{c}_1},
\end{align}
where $\hat{\mathscr{C}}^\dagger_{\!\!\boldsymbol{X}^{[L]}} \coloneqq (\hat{\mathscr{C}}_{\!\boldsymbol{X}^{[L]}} + \theta_{\!N} \mathbb{I})^{-1}$ is the Tikhonov-type pseudo-inverse of $\hat{\mathscr{C}}_{\!\boldsymbol{X}^{[L]}},$ with identity map $\mathbb{I}$ and positive tuning parameters with $\theta_{\!N} \to 0,$ and  $K = K_N \to \infty.$ The eigenfunctions $\hat{c}_1, \dots, \hat{c}_K$ correspond to the $K$ largest eigenvalues $\hat{\lambda}_1 \geq \dots \geq \hat{\lambda}_K$ of $\hat{\mathscr{C}}_{\!\boldsymbol{X}^{[L]}}$. To analyze the error between $\hat{\Psi}_{\!L}$ and $\Psi_{\!L}$, we decompose $\Psi_{\!L}$ into a finite-dimensional component and a remainder term, controlled via a Sobolev condition similar to that in  \cite{hall:meister:2007}.  Let $c_i = c_i(L)$ and $d_j$ be the eigenfunctions of $\mathscr{C}_{\!\boldsymbol{X}^{[L]}} \colon \mathcal{H}^L \to \mathcal{H}^L$ and $\mathscr{C}_{\!\boldsymbol{X}} \colon \mathcal{H} \to \mathcal{H}$, respectively. By Assumption \ref{As: Inv process}, $\Psi_{\!L} \colon \mathcal{H}^L \to \mathcal{H}$ is H-S, and admits the decomposition  
$$
\Psi_{\!L}=\sum_{i=1}^\infty\sum_{j=1}^\infty\,\langle\Psi_{\!L},c_i\otimes d_j\rangle_{\mathcal{S}}(c_i\otimes d_j).
$$

\begin{assumption}\label{As: Sob cond inv ops} For some $\beta >0,$ the limit $S_{\Psi}(\beta) \coloneqq \lim_{L \rightarrow\infty}S_{\Psi_{\!L}}(\beta) < \infty$ exists, where for each $L=L_N\in\mathbb{N}$,
  \begin{gather*}
   S_{\Psi_{\!L}}\!(\beta) 
   \coloneqq \sum_{i=1}^\infty\sum_{j=1}^\infty\,\langle\Psi_{\!L}, c_i\otimes d_j\rangle^2_\mathcal{S} \,(1 + i^{2\beta}\! + \hspace{-0.05em}j^{2\beta}) <
	\infty.
  \end{gather*}
\end{assumption}
\noindent Our results bellow allow for repeated eigenvalues $\lambda_j=\lambda_j(L)$ of $\mathscr{C}_{\!\boldsymbol{X}^{[L]}},$ provided the dimension of each associated eigenspace $\mathcal{V}_j(L)$ of $\mathscr{C}_{\!\boldsymbol{X}^{[L]}}$ is bounded in the following sense.

\begin{assumption}\label{As: Dim eigenspaces unif bounded inv op} There exists a constant $d\in\mathbb{N}$ so that  $\sup_j\mathrm{dim}(\mathcal{V}_j(L)) \leq dL.$ 
\end{assumption}

\noindent In order to state our main result, we define for each $L$ a real-valued sequence $\Lambda_\ell = \Lambda_\ell (L)$, $\ell \in \mathbb{N}$, related to the reciprocal eigengaps of  $\mathscr{C}_{\!\boldsymbol{X}^{[L]}}$.  Let  
$$
\ell_b = \inf \left\{  j \colon  j > \ell, \; \lambda_j(L) < \lambda_\ell(L) \right\}.
$$
Define 
\begin{align}\label{eq:Def Lambda_K}
\Lambda_\ell \coloneqq \frac{1}{\lambda_\ell - \lambda_{\ell_b}}.
\end{align}
In other words, $\Lambda_\ell$ is the reciprocal eigengap between the eigenspace associated to $c_\ell$, and the one after it. Another assumption we need is that the empirical eigenvalues do not occur repeatedly, so that the empirical reciprocal eigengaps are well-defined.

\begin{assumption}\label{As:No rep eigenvalues main Theo} For all $N$ sufficiently large, with probability one the eigenvalues $\hat{\lambda}_j=\hat{\lambda}_j(L)$ of $\hat{\mathscr{C}}_{\!\boldsymbol{X}^{[L]}}$ satisfy $\hat{\lambda}_j\neq \hat{\lambda}_{j+1}$ for all $j\in\{1,...,K+dL\}$ and $L,$ with $d$ from Assumption \ref{As: Dim eigenspaces unif bounded inv op}.
\end{assumption}

\begin{theorem}\label{Theo: inv. op. - infinite-dim} Let Assumptions \ref{As: Features of (eps_k)}--\ref{As:No rep eigenvalues main Theo} hold. Further, let $K=K_{\!N} \to \infty$ and $\theta_{\!N} \to 0$ be sequences so that for $L=L_{\!N}=\mathrm{o}(N^{1/3})$ from Assumption \ref{As: cons rate for lagged cov and cross-cov ops} holds $L=\mathrm{O}(K),$ and 
\begin{align*}
  K^{\beta+1/2}\Lambda^2_KL^{2} = \mathrm{O}(N^{1/2}), \quad \sum_{\ell>L}\|\psi_\ell\|_\mathcal{L} = \mathrm{O}(\lambda_K\Lambda^2_KK^{1/2}L^{3/2}N^{-1/2}), \quad \theta_{\!N} = \mathrm O(\lambda_KK^{-\beta}).
\end{align*}
Then, 
\begin{align*}
  \|\hat{\Psi}_{\!L} - \Psi_{\!L}\|_\mathcal{S}  = \,\mathrm O_{\Prob}(K^{-\beta}).
\end{align*}
\end{theorem}

\begin{remark}\label{Remark on operators in inverted representation}
\begin{itemize}
    \item[\textnormal{(a)}] Theorem \ref{Theo: inv. op. - infinite-dim} advances operator estimation for invertible processes by providing explicit rates for the H-S norms of the estimation errors for the entire operators in any separable Hilbert space, allowing multi-dimensional eigenspaces. Previous works, mainly in $\mathcal{H} = L^2[0,1]$, often focused on finite-dimensional approximations, omitted rates in infinite dimensions, or assumed one-dimensional eigenspaces. For related results in f(G)ARCH models, see \cite{Hoermannetal2013, Aueetal2017, Ceroveckietal2019, Kuehnert2020}; and for invertible processes, see \cite{AueKlepsch2017, klepsch2017innovations}.
\item[\textnormal{(b)}] The rate in Theorem \ref{Theo: inv. op. - infinite-dim} is primarily governed by the decay  of the eigenvalues $\lambda_K = \lambda_K(L)$ of $\mathscr{C}_{\! \boldsymbol{X}^{[L]}}$, the smoothness of $\psi_j$ (Assumption \ref{As: Sob cond inv ops}), and, to a lesser extent, the decay of $\|\psi_j\|_\mathcal{L}$. Estimation is most effective when the $\psi_j$ are smooth, and  the norms $\|\psi_j\|_\mathcal{L}$ decay rapidly.
Notably, if the eigenvalues of $\psi_j$ decay geometrically, the parametric rate $N^{-1/2}$ can be almost reached (see Examples \ref{Ex: Explicit process w/ rates fAR1} and \ref{Ex: Explicit process w/ rates fARinfty}).

    \item[\textnormal{(c)}] Assumption \ref{As: Sob cond inv ops} ensures square-summability of $\|\psi_j\|_\mathcal{S}$ via Parseval’s identity. Assumption \ref{As: Dim eigenspaces unif bounded inv op} imposes minimal restrictions and even accomodates the case that $\boldsymbol{X} = (X_k)$ is i.i.d. \frenchspacing if the eigenspaces of the $\mathscr{C}_{\!\boldsymbol{X}}'$ have dimension up to $d,$ since $\mathscr{C}_{\! \boldsymbol{X}^{[L]}} = \mathrm{diag}(\mathscr{C}_{\!\boldsymbol{X}}, \dots, \mathscr{C}_{\!\boldsymbol{X}}).$ For conditions on the eigengaps $f(j) = \lambda_j - \lambda_{j-1}$, see, e.g.,  \cite{hall:horo:2007}, equation 3.2.
\end{itemize}
\end{remark}

We conclude this section by illustrating the interplay between the conditions in Theorem \ref{Theo: inv. op. - infinite-dim}. 

\begin{example}\label{Ex: Explicit process w/ rates fAR1} 
Consider the invertible process $(X_k)_{k\in\mathbb{Z}} \subset \mathcal{H}$ given by
\begin{equation*}
  X_k =  \sum_{j=1}^{\infty} \psi^j(X_{k-j}) + \varepsilon_k,\quad k\in\mathbb{Z},
\end{equation*}
where $(\varepsilon_k)_{k \in \mathbb{Z}} \subset L^4_\mathcal{H}$ is an i.i.d. process with injective covariance operator $\mathscr{C}_{\boldsymbol{\varepsilon}}$, having eigenvalues $a_\ell$ and eigenfunctions $e_\ell$. Suppose the operator $\psi$ admits the spectral representation
\begin{equation}\label{svd of psi}
  \psi = \sum^\infty_{\ell=1} p_\ell(e_\ell\otimes e_\ell),    
\end{equation}
with $1/2>p_1\geq p_2 \geq \cdots >0,$ and $p_\ell\to0$. Since $\|\psi\|_\mathcal{L}=p_1<1/2$, iterative application of the invertible representation yields the linear representation
\begin{equation}\label{Lin repr of inv proc in Ex}
    X_k = \sum_{j=1}^{\infty}\, 2^{j-1}\psi^j(\varepsilon_{k-j}) + \varepsilon_k,\quad k\in\mathbb{Z}.
\end{equation}
The series converges a.s. \frenchspacing and in $L^2_\mathcal{H}$, ensuring weak stationarity \citep[][Section 7]{Bosq2000}. As we assume $(\varepsilon_k)$ is i.i.d., $(X_k)$ is $L^4$-$m$-approximable \citep[][Proposition 2.1]{hormann:kokoszka:2010}, satisfying Assumption \ref{As: cons rate for lagged cov and cross-cov ops} \citep[][Proposition 3.2 and Corollary 3.1]{Kuehnert2025}.
According to \eqref{Lin repr of inv proc in Ex} and since $(\varepsilon_k)$ is i.i.d., the lag-$h$-covariance operators have the form
\begin{equation*}
  \mathscr{C}^h_{\!\boldsymbol{X}} = \mathbb{E}(X_0\otimes X_h) = \bigg[\,2^{|h|-2}\psi^{|h|}\sum^\infty_{i=1} \psi^i + 2^{|h|-1}\psi^{|h|}\cdot\mathds{1}_{\mathbb{Z}\backslash\{0\}}(h) + \mathbb{I}\cdot\mathds{1}_{\{0\}}(h)\bigg]\mathscr{C}_{\boldsymbol{\varepsilon}}, \quad h\in\mathbb{Z}.
\end{equation*} 
Thus, $\mathscr{C}_{\!\boldsymbol{X}}=\mathscr{C}^0_{\!\boldsymbol{X}}$ and $\mathscr{C}^h_{\!\boldsymbol{X}}$ share eigenfunctions $e_\ell$ with $\mathscr{C}_{\!\boldsymbol{\varepsilon}}$, with eigenvalues
\begin{equation*}
  \mu^0_\ell = \frac{a_\ell}{2-4p_\ell}(2-3p_\ell), \quad \mu^h_\ell = \frac{a_\ell}{2-4p_\ell} (1-p_\ell)(2p_\ell)^{|h|}, \quad h\neq 0.
\end{equation*}

\noindent The operator-valued matrix $\mathscr{C}_{\! \boldsymbol{X}^{[L]}} = (\mathscr{C}^{j-i}_{\! \boldsymbol{X}})^L_{i,j}$ has eigenfunctions
\begin{equation*}
    c_{\ell }= e_{\lfloor (\ell-1)/L \rfloor +1}v_{\ell, ((\ell-1)\,\mathrm{mod}\,L) +1},
\end{equation*}
with eigenvectors $v_{\ell, 1}, \dots, v_{\ell, L}$ corresponding to eigenvalues $t_{\ell,1}\geq \dots \geq t_{\ell,L}$ of $T_\ell = (\mu^{|j-i|}_\ell)^L_{i,j}.$ By Gershgorin's theorem,
\begin{equation*}
    \big|t_{\ell, m} - \mu^0_\ell\big| \leq  \max_{1\leq i \leq L} \sum^L_{j=1, j\neq i}\mu^{|j-i|}_\ell = \sum^{L-1}_{j=1}\mu^{|j-\lfloor(L+1)/2\rfloor|}_\ell,\quad m=1, \dots, L.
\end{equation*}
Since this sum converges as $L\to\infty$ and is of order $\mathrm{O}(a_\ell)$ as $\ell\to\infty$ (as is $\mu^0_\ell$), the eigenvalues of $\mathscr{C}_{\! \boldsymbol{X}^{[L]}}$ satisfy $\lambda_\ell(L) \asymp a_{\lfloor \ell/L \rfloor},$ ensuring Assumption \ref{As: Dim eigenspaces unif bounded inv op}. For Assumption \ref{As: Sob cond inv ops}, setting $d_j=e_j$ and $\Psi_{\!L}=(\psi^1, \dots, \psi^L)^\top$ yields $\langle\Psi_{\!L}, c_i \otimes d_j\rangle_\mathcal{S} \leq \delta_{ij} \sum_{m\leq L}p^m_i,$ where $\delta_{ij}$ is the Kronecker delta. Hence,  for $\beta<b-1/2$ if $p_\ell \asymp \ell^{-b}$ for some $b>1$, or any $\beta>0$ if $p_\ell \asymp q^\ell$ for some $q\in(0,1),$ it holds that
\begin{equation*}
    \lim_{L\to\infty}S_{\Psi_{\!L}}(\beta)=  \lim_{L\to\infty}\sum_{i=1}^{\infty} \sum_{j=1}^{\infty}\,\langle\Psi_{\!L}, c_i \otimes d_j\rangle^2_\mathcal{S} (1 + i^{2\beta} + j^{2\beta}) \leq \sum_{i=1}^{\infty}\,\frac{p^2_i(1 + 2i^{2\beta})}{(1-p_i)^2}<\infty.
\end{equation*}

It remains to determine the sequences in Theorem \ref{Theo: inv. op. - infinite-dim}. With $\theta_{\!N}$ chosen freely within the given constraints, we know $\sum_{\ell>L} \|\psi_\ell\|_\mathcal{L}\propto \|\psi\|^L_\mathcal{L} = \mathrm{o}(N^{-1/2}),$ whenever $L$ is a positive power of $N$. Valid choices of $K=K_{\!N}$ with  $K^{\beta+1/2}\Lambda^2_KL^2  \asymp N^{1/2}$ is for given $\beta$ and $L$ determined by the decay of the eigenvalues $\lambda_j.$ E.g., for  $a_\ell\asymp \ell^{-a}$ for some $a>1,$ it holds $\lambda_K \asymp (KL^{-1})^{-a}$ and thus $\Lambda_K \asymp (KL^{-1})^{a+1}.$ Hence, with $L=L_{\!N}\asymp K^c$ for small $c\in(0,1)$ and $K\asymp N^{1/(5+4a(1-c)+2\beta)},$ Theorem \ref{Theo: inv. op. - infinite-dim} yields
\begin{equation*}
    \|\hat{\Psi}_{\!L} - \Psi_{\!L}\|_\mathcal{S} 
    = \mathrm{O}_{\Prob}\Big(N^{-\frac{\beta}{5+4a(1-c)+2\beta}}\Big).
\end{equation*}
For e.g., $a_\ell\asymp \ell^{-2}$ \citep[the eigenvalue decay rate for Brownian motion errors, see][]{JaimezBonnet1987}, the rates are, for instance, close to $N^{-1/8}$ for $p_\ell \asymp \ell^{-13/6}$, $N^{-1/4}$ for $p_\ell \asymp \ell^{-13/2}$, and approximately the parametric rate $N^{-1/2}$ for $p_\ell \asymp \rho^\ell$, $\rho \in (0,1)$.
\end{example}

\section{Consequences for causal linear processes}\label{Section: Linear process representation}

We now turn to the estimation of operators in the linear representation of invertible, causal linear processes. Let $(X_k)_{k\in\mathbb{Z}} \subset \mathcal{H}$ follow a causal linear process (see \eqref{Definition LP} for general cases),  
\begin{align}\label{identity of linear process}
  X_k =  \sum^\infty_{j=0} \phi_j(\varepsilon_{k-j}), \quad k\in\mathbb{Z},
\end{align}
where $(\varepsilon_k)_{k\in\mathbb{Z}} \subset \mathcal{H}$ is a WN, and $\phi_j \in \mathcal{L}_{\mathcal{H}}$ are operators with $\phi_0 \coloneqq \mathbb{I}$ and $\sum_{j=1}^{\infty} \|\phi_j\|_\mathcal{L} < \infty.$ The process is \emph{invertible} if it admits a representation as in \eqref{identity of invertible process}. For a comprehensive treatment of linear processes in function spaces, see \cite{Bosq2000}, particularly Section 7 for existence and invertibility results. We now focus on estimating the operators $\phi_j,$ which is crucial for prediction and confidence set construction in $\mathcal{H}$-valued models. The relationship between estimates from an inverted representation of $(X_k)$ and its causal form follows from the next proposition.  

\begin{proposition}\label{Prop: Unique operators in LP and inv process} Let Assumption \ref{As: Features of (eps_k)} hold. Then the following is true.
  \begin{itemize}
    \item[\textnormal{(a)}] If for some $i \in \mathbb{N}$ it holds that $\sum^\infty_{\ell=0}A_\ell(\varepsilon_{i-\ell}) = \sum^\infty_{\ell=0}B_\ell(\varepsilon_{i-\ell}),
    $ then $A_\ell = B_\ell$ for all $\ell\in \mathbb{N}.$
    \item[\textnormal{(b)}] If $(X_k)$ is a causal linear process and if for some $i\in\mathbb{N}$, $\sum^\infty_{\ell=0}A_\ell(X_{i-\ell}) = \sum^\infty_{\ell=0}B_\ell(X_{i-\ell}),$ then $A_\ell = B_\ell$ for all $\ell \in \{0,1,2, \dots\}.$
  \end{itemize}
\end{proposition}

For an invertible linear process $(X_k)$ in Assumption \ref{As: Inv process} such that

$$X_k = \sum^\infty_{i=1}\psi_i(X_{k-i}) + \varepsilon_k,\quad \mbox{ and } \quad X_{k-i}=\sum^\infty_{j=1} \phi_j(\varepsilon_{k-i-j}) + \varepsilon_{k-i},$$
for all $i,k \in \mathbb{Z}$, we have due to $\phi_0=\mathbb{I}$ that
\begin{align*}
X_k &= \sum^\infty_{i=1}\,\Big[\sum^i_{j=1}\psi_j\phi_{i-j}\Big](\varepsilon_{k-i}) + \varepsilon_k,\quad k\in\mathbb{Z}.
\end{align*}
\noindent This representation, the linear representation \eqref{identity of linear process}, and Proposition \ref{Prop: Unique operators in LP and inv process} (a) yield
\begin{align*}
  \phi_i = \sum^i_{j=1}\psi_j\phi_{i-j}, \quad i \in \mathbb{N}.
\end{align*}

\noindent Hence, with $\hat{\psi}_{j} =\hat{\Psi}^{(j)}_{\!L}$ denoting the $j$'th component of the estimator $
\hat{\Psi}_{\!L}$ in \eqref{Operator estimate for big Psi_L} for $\Psi_{\!L}=(\psi_1, \dots, \psi_L),$ a reasonable estimate for $\phi_i$ is $\hat{\phi}_i,$ defined iteratively with $\hat{\phi}_0\coloneqq \mathbb{I},$  and
\begin{align}\label{Definition of estimators for phi_i}
  \hat{\phi}_i\coloneqq \sum^i_{j=1}\,\hat{\psi}_{j}\hat{\phi}_{i-j}, \quad i \in \mathbb{N}.
\end{align}

The consistency result for the complete operators in the inverted representation carries over to those in the linear representation.

\begin{proposition}\label{Proposition linear operator estimation - infinite-dim} Let the conditions of Theorem \ref{Theo: inv. op. - infinite-dim} hold, and let $\phi_i\in\mathcal{S_H}$ for all $i.$ Then, 
  \begin{align*}
    \|\hat{\phi}_i - \phi_i\|_\mathcal{S} = \mathrm{O}_{\Prob}\big(\|\hat{\Psi}_{\!L} - \Psi_{\!L}\|_\mathcal{S}\big)  = \mathrm O_{\Prob}(K^{-\beta}),\quad i\in\mathbb{N}.
  \end{align*}
\end{proposition}

\section{Consequences for fARMA processes}\label{Section: Implications to ARMA}

The asymptotic results in Section \ref{Section: Estimation in the invertible and linear representation} extend to the estimation of fAR, fMA, and fARMA parameters. These models are extensively studied in \cite{Bosq2000}, particularly the fAR case, while \cite{Spangenberg2013} provides general conditions for causal and invertible representations in Banach spaces. Estimation results for the fAR(1) operator are available in \cite{Bosq2000} and \cite{CaponeraPanaretos2022}. Additionally, \cite{TurbillonETAL2008} estimated the fMA(1) operator assuming commutativity with the error covariance operator, and \cite{AueKlepsch2017, KlepschKlueppelbergWei2017} studied prediction for fARMA$(p,q)$ processes. \cite{kuenzer:2024} derived consistent fMA(1) and fARMA(1,1) operator estimates using an adaptation of the Hannan–Rissanen method \citep{hannan:1982}. By using operator estimates from both the linear and inverted representations, we derive consistent fAR, fMA, and fARMA operator estimates with explicit rates.

\subsection{fAR processes}
A centered process $(X_k)_{k\in\mathbb{Z}}\subset  \mathcal{H}$ is a \emph{functional autoregressive process} of order $p\in\mathbb{N}$ (fAR$(p)$) if
\begin{align*}
  X_k= \sum^{p}_{i=1}\,\alpha_i(X_{k-i}) + \varepsilon_k, \quad k \in \mathbb{Z},
\end{align*}
where $(\varepsilon_k)_{k\in\mathbb{Z}}\subset\mathcal{H}$ is a WN, and $\alpha_i\in \mathcal{L_H}$ are operators with $\alpha_p\neq 0.$   With $\Psi_{\!p}$ in \eqref{Definition of operator-valued 1xL-matrix}, it holds 
$\Psi_{\!p} = (\alpha_1\,\cdots\,\alpha_p),$ and the exact Yule--Walker equation $\mathscr{C}^1_{\!\boldsymbol{X}^{[p]}, \boldsymbol{X}} = \Psi_{\!p}\mathscr{C}_{\!\boldsymbol{X}^{[p]}}$ is fulfilled. Then, with
$$
\hat{\alpha}_i \coloneqq \hat{\Psi}^{(i)}_{\!p}, \quad i=1, \dots, p,
$$
which are the $i$'th components of $\hat{\Psi}_{\!p}$ in \eqref{Operator estimate for big Psi_L}, we obtain the following corollary of Theorem \ref{Theo: inv. op. - infinite-dim}.

\begin{corollary}\label{Corollary AR(p) - finite-dim} Let the conditions of Theorem \ref{Theo: inv. op. - infinite-dim} with $L=p$ hold, and let $\alpha_i \in \mathcal{S_H}$ for all $i.$ Then, 
\begin{gather*}
    \max_{1\leq i \leq p}\|\hat{\alpha}_i - \alpha_i\|_\mathcal{S}  = \mathrm O_{\Prob}(K^{-\beta}). 
\end{gather*}
\end{corollary}

\begin{example}\label{Ex: Explicit process w/ rates fARinfty} Let $(X_k)_{k\in\mathbb{Z}} \subset \mathcal{H}$ be the fAR(1) process  
\begin{equation*}
    X_k = \alpha_1(X_{k-1}) + \varepsilon_k, \quad k \in \mathbb{Z},
\end{equation*}
where $\boldsymbol{\varepsilon}=(\varepsilon_k)_{k \in \mathbb{Z}} \subset L^4_\mathcal{H}$ is i.i.d. with injective covariance operator $\mathscr{C}_{\boldsymbol{\varepsilon}}$ having eigenvalues $a_\ell$. The operator $\alpha_1$ has the same form as $\psi$ in Example \ref{Ex: Explicit process w/ rates fAR1}. Since $\|\alpha_1\|_\mathcal{L} < 1$, Theorems 3.1--3.2 in \cite{Bosq2000} ensure a unique weakly stationary solution, where the covariance operator has the form
\begin{gather*}
    \mathscr{C}_{\!\boldsymbol{X}} = \mathscr{C}_{\boldsymbol{\varepsilon}} \sum_{i=0}^{\infty} \alpha_1^{2i}, 
\end{gather*}  
with eigenvalues $\lambda_\ell \propto a_\ell.$ To derive the asymptotic rate, we follow Example \ref{Ex: Explicit process w/ rates fARinfty}, setting $L=p=1.$ Assumption \ref{As: Sob cond inv ops} is fulfilled if the eigenvalues $p_i$ of $\alpha_1$ fulfill $\sum_{i=1}^{\infty} i^{2\beta}p_i^2< \infty.$ Then, if $a_\ell \asymp \ell^{-a}$ for some $a > 1$, e.g., Corollary \ref{Corollary AR(p) - finite-dim} ensures that $K$ can be chosen so that  
\begin{equation*}
    \|\hat{\alpha}_1 - \alpha_1\|_\mathcal{S} 
    = \mathrm{O}_\mathbb{P}\Big( N^{-\frac{\beta}{5+4a+2\beta}} \Big).
\end{equation*}  
\end{example}

\subsection{fMA processes}
A centered process $(X_k)_{k\in\mathbb{Z}} \subset \mathcal{H}$ is a \emph{functional moving-average process} of order $q \in \mathbb{N}$ (fMA$(q)$) if  
\begin{align*}
  X_k=  \sum^q_{j=1}\,\beta_j(\varepsilon_{k-j}) + \varepsilon_k, \quad k \in \mathbb{Z},
\end{align*}
where $(\varepsilon_k)_{k\in\mathbb{Z}}\subset\mathcal{H}$ is a WN, and $\beta_j\in \mathcal{L_H}$ are operators with $\beta_q\neq 0.$ Morever, let 
$$
  \hat{\beta}_j = \hat{\phi}_j, \quad j=1, \dots, q,
$$
with $\hat{\phi}_j$ in \eqref{Definition of estimators for phi_i}. Based on the previous considerations, we immediately obtain the following result. 
\begin{corollary}\label{Corollary from linear operators} Let the conditions of Theorem \ref{Theo: inv. op. - infinite-dim} hold, and let $\beta_j \in \mathcal{S_H}$ for all $j.$ Then, it holds  \begin{align*}
    & \max_{1\leq j \leq q}\|\hat{\beta}_j - \beta_j\|_\mathcal{S} = \mathrm O_{\Prob}(K^{-\beta}).
  \end{align*}
\end{corollary}

\subsection{fARMA processes}\label{Subsec: fARMA}

We finally turn to the estimation of fARMA processes. A centered process $(X_k)_{k\in\mathbb{Z}}\subset \mathcal{H}$ is called a   \emph{functional autoregressive moving-average process} of order $p,q\in\mathbb{N}$ (fARMA$(p, q)$) if
\begin{align}\label{functional ARMA definition}
  X_k= \sum^{p}_{i=1}\,\alpha_i(X_{k-i})+\sum^{q}_{j=1}\,\beta_j(\varepsilon_{k-j}) + \varepsilon_k , \quad k \in \mathbb{Z},
\end{align}
where $(\varepsilon_k)_{k\in\mathbb{Z}}\subset\mathcal{H}$ is a WN, and $\alpha_i, \beta_j \in \mathcal{L_H}$ are operators with $\alpha_p\neq 0, \beta_q \neq 0.$ We combine this fARMA$(p,q)$ representation with its inverted form \eqref{identity of invertible process}. As with multivariate ARMA processes, the operators $\alpha_i$ and $\beta_j$ are generally not identifiable. Sufficient conditions for identifiability are provided in \cite{kuenzer:2024} (equation (10) and Propositions 1–3). We assume that the fARMA models considered are identifiable and admit invertible, centered, linear process solutions.

To describe the main idea, we first discuss estimation for fARMA$(1,1)$ processes of the form
\begin{align}\label{fARMA(1,1) definition}
  X_k= \alpha_1(X_{k-1})+\,\beta_1(\varepsilon_{k-1}) + \varepsilon_k, \quad k \in \mathbb{Z}.
\end{align}
\noindent From this identity, Proposition \ref{Prop: Unique operators in LP and inv process} and $\varepsilon_{k} = X_k - \sum^\infty_{i=1}\psi_i(X_{k-i})$ it follows with $\psi_0=\mathbb{I},$
\begin{gather}\label{equation of invertible operators in ARMA(1,1)}
\psi_i =
  \begin{cases}
    \,\alpha_1 + \beta_1, &\textnormal{if}\;\,i=1,\\
	\,-\beta_1\psi_{i-1}, &\textnormal{if}\;\,i > 1.\\
  \end{cases}
\end{gather}

Although $\hat{\Psi}_{\!L}$ in \eqref{Operator estimate for big Psi_L} provides consistent asymptotic estimates for all $\psi_i$, consistent estimates for the fARMA$(1,1)$ operators cannot be directly obtained, as $\psi_1 \in \mathcal{S_H}$ does not possess a bounded inverse. However, using $\hat{\psi}_i = \hat{\Psi}^{(i)}_{\!L}$, equation \eqref{equation of invertible operators in ARMA(1,1)} leads to the following estimators for $\alpha_1$ and $\beta_1$:
\begin{gather}\label{Estimator for alpha_1 and beta_1}
  \hat{\alpha}_1 \coloneqq \hat{\psi}_1 - \hat{\beta}_1, \qquad
  \hat{\beta}_1 \coloneqq -\hat{\psi}_2 \hat{\psi}_1^\dagger \coprod^{\hat{f}_M}_{\hat{f}_1},
\end{gather}
where $\hat{\psi}_1^\dagger \coloneqq (\hat{\psi}_1 \hat{\psi}_1^\ast + \gamma_N \mathbb{I})^{-1},$ and $(M_N) \subset \mathbb{N}$ and $(\gamma_N) \subset (0, \infty)$ are sequences with $M_N \to \infty$ and $\gamma_N \to 0$. Additionally, $\hat{f}_1, \dots, \hat{f}_M$ are the eigenfunctions corresponding to the eigenvalues $\hat{\rho}_1 \geq \dots \geq \hat{\rho}_M \geq 0$ of $\hat{\psi}_1 \hat{\psi}_1^\ast \colon \mathcal{H} \to \mathcal{H}$. The sequence $(\rho_j, f_j)$ represents the eigenpairs of the operator $\psi_1 \psi_1^\ast \colon \mathcal{H} \to \mathcal{H}$. The estimates in \eqref{Estimator for alpha_1 and beta_1} are consistent under the following assumptions, where the first assumption ensures identifiability of $\beta_1$ from $\psi_2 = -\beta_1 \psi_1.$

\begin{assumption}\label{As: psi_1psiadj_1 injective} The image of the operator $\psi_1\colon\mathcal{H}\rightarrow \mathcal{H}$ lies dense.
\end{assumption}

\begin{assumption}\label{As: Sobolev cond ARMA(1,1) - alpha_1beta_1} The operators $\alpha_1, \beta_1\in\mathcal{S_H}$ are H-S, and  for some $\gamma >0,$ it holds 
  \begin{gather*}
  \sum_{i=1}^\infty\sum_{j=1}^\infty\,\langle\beta_1, f_i\otimes f_j\rangle^2_\mathcal{S}(1 + i^{2\gamma}\! + \hspace{-0.05em}j^{2\gamma}) < \infty.
  \end{gather*}
\end{assumption}

\begin{assumption}\label{As: Dim and eigvals fARMA 11} There exists $d \in \mathbb{N}$ such that:
\begin{itemize}
    \item[\textnormal{(a)}] The dimensions of all eigenspaces of $\psi_1\psi^\ast_1$ are bounded above by $d;$
    \item[\textnormal{(b)}] For all $N$ sufficiently large, with probability one the eigenvalues of $\hat{\psi}_1\hat{\psi}^\ast_1$ satisfy $\hat{\rho}_j\neq \hat{\rho}_{j+1}$ for all $j\in\{1,...,M+d\}.$
\end{itemize}
\end{assumption}

\begin{theorem}\label{Theorem for ARMA(1,1) operators, infinite-dim} Suppose $(X_k)$ is  the fARMA$(1,1)$ process in \eqref{fARMA(1,1) definition}. Further, let the conditions of Theorem \ref{Theo: inv. op. - infinite-dim} and Assumptions \ref{As: psi_1psiadj_1 injective}--\ref{As: Dim and eigvals fARMA 11} hold, and let $M=M_{\!N}\rightarrow\infty$ with $
M^{\gamma+1/2}P^2_{\!M} = \mathrm{O}(K^\beta)$ and $\gamma_N = \mathrm O(\rho_MM^{-\gamma}),$ 
with $P_M\coloneqq (\rho_M - \rho_{M_b})^{-1}$ as in \eqref{eq:Def Lambda_K}. Then, 
\begin{align*}
  \max\!\big\{\|\hat{\alpha}_1 - \alpha_1\|_\mathcal{S}, \|\hat{\beta}_1 - \beta_1\|_\mathcal{S}\big\} = \mathrm O_{\Prob}(M^{-\gamma}).
  \end{align*}
\end{theorem}

The procedure for arbitrary $p,q$ is similar as for $p=q=1.$ Similar to the transformations that led to \eqref{equation of invertible operators in ARMA(1,1)}, we obtain for any $p$ and $q,$ with $\alpha_i = \beta_j = 0$ for $i>p, j>q,$
\begin{gather}
  \sum^\infty_{i=1}\,\psi_i(X_{k-i})\, = \sum^{\max(p,q)}_{i=1}\Big[\alpha_i+\beta_i - \sum^{i-1}_{j=1}\,\beta_j\psi_{i-j}\Big](X_{k-i})\, - \sum^{\infty}_{i=\max(p,q)+1}\Big[\sum^q_{j=1}\,\beta_j\psi_{i-j}\Big](X_{k-i})\notag.
\end{gather}
Comparing the left- and right-hand sides of the above and applying Proposition \ref{Prop: Unique operators in LP and inv process} (b), we see that
\begin{gather}
 \psi_i =
  \begin{cases}
    \,\alpha_i+\beta_i - \sum^{i-1}_{j=1}\,\beta_j\psi_{i-j}, &\textnormal{if}\;\,1\leq i \leq \max(p,q),\\
	\,-\sum^q_{j=1}\beta_j\psi_{i-j}, &\textnormal{if}\;\,i > \max(p,q).\\
  \end{cases}\label{equation of invertible operators in ARMA(p,q)}
\end{gather}
We proceed to obtain estimates for $\beta_i$ for $i>\max(p,q),$ and then for $\alpha_i$. Since \eqref{equation of invertible operators in ARMA(p,q)} involves sums of compositions of the $\beta_j$ with $\psi_{i-j}$ for $i>\max(p,q),$ we cannot immediately retrieve estimates for $\beta_j$ by applying Tikhonov regularized versions of the estimates $\psi_{i-j}$. For $i=p+q,$ with
  \begin{gather}
 B_q \coloneqq (\beta_1\,\cdots\,\beta_q) \in \,\mathcal{L}_{\mathcal{H}^q\!,\mathcal{H}},\notag\\
  \Psi'_{[i]}\coloneqq (\psi_{p+q+i-1}\,\cdots\, \psi_{p+i}) \in \,\mathcal{L}_{\mathcal{H}^q\!,\mathcal{H}},\quad i \geq 0,\label{definition of vector of psis for ARMA(p,q)}
\end{gather}
the identity \eqref{equation of invertible operators in ARMA(p,q)} becomes  $\psi_{p+q} =  -B_q \Psi'^{\,\top}_{[0]}.$ Identifiability of $B_q$ and thus of $\beta_1, \dots, \beta_q$ is given if the image of $\Psi'^{\,\top}_{[0]}\in \,\mathcal{L}_{\mathcal{H}, \mathcal{H}^q}$ would lie  dense. This condition is never met. Nonetheless, identifiability of $B_q$ can be established through the following relationship:
\begin{align}\label{equation of invertible operators in ARMA(p,q), vector representation}
   \Psi'_{[q]}  =  -B_q\,\boldsymbol{\prod}\;,
\end{align}
where $\boldsymbol{\prod}\in \mathcal{L}_{\mathcal{H}^q}$ is the operator-valued matrix defined by
\begin{align}
  \boldsymbol{\prod}\,\coloneqq \begin{bmatrix}
\,\Psi'_{\![q-1]}\,\\
\,\Psi'_{\![q-2]}\,\\
\,\vdots\,\\
\,\Psi'_{\![0]}\,\\
\end{bmatrix} = \begin{bmatrix}
\,\psi_{p+2q-2} & \psi_{p+2q-3} & \cdots & \psi_{p+q-1}\,\\
\,\psi_{p+2q-3} & \psi_{p+2q-4} & \cdots & \psi_{p+q-2}\,\\
\,\vdots   & \vdots      & \cdots     & \vdots\,\\
\,\psi_{p+q-1} & \psi_{p+q-2} & \cdots & \psi_p\\
\end{bmatrix}.\label{definition of matrices of psis for ARMA(p,q)}
\end{align}
For $\Psi'_{[i]}$ and $\boldsymbol{\prod},$ we use the estimates $\hat{\Psi}'_{[i]}$ and $\boldsymbol{\hat{\prod}},$ respectively, which are defined by
\begin{gather}
    \hat{\Psi}'_{[i]}\coloneqq (\hat{\psi}_{p+q+i-1}\,\cdots\, \hat{\psi}_{p+i}) \in \,\mathcal{L}_{\mathcal{H}^q\!,\mathcal{H}},\quad i \geq 0, \mbox{ and }\label{definition of estimate for vector of psis for ARMA(p,q)}\\
  \boldsymbol{\hat{\prod}}\,\coloneqq \begin{bmatrix}
\,\hat{\Psi}'_{\![q-1]}\,\\
\,\hat{\Psi}'_{\![q-2]}\,\\
\,\vdots\,\\
\,\hat{\Psi}'_{\![0]}\,\\
\end{bmatrix} = \begin{bmatrix}
\,\hat{\psi}_{p+2q-2} & \hat{\psi}_{p+2q-3} & \cdots & \hat{\psi}_{p+q-1}\,\\
\,\hat{\psi}_{p+2q-3} & \hat{\psi}_{p+2q-4} & \cdots & \hat{\psi}_{p+q-2}\,\\
\,\vdots   & \vdots      & \cdots     & \vdots\,\\
\,\hat{\psi}_{p+q-1} & \hat{\psi}_{p+q-2} & \cdots & \hat{\psi}_p\\
\end{bmatrix}.
\label{definition of estimate for matrices of psis for ARMA(p,q)}
\end{gather}
Due to Lemma \ref{Lemma: Inequalities norm of operator matrices} (b), the operator-valued vectors $\Psi'_{[i]}, \hat{\Psi}'_{[i]}$ and matrices $\prod, \hat{\prod}$ are under Assumption \ref{As: Sob cond inv ops}  elements of $\mathcal{S}_{\mathcal{H}^q\!,\mathcal{H}}$ and $\mathcal{S}_{\mathcal{H}^q},$ respectively, thus $\prod\prod^\ast, \hat{\prod}\,\hat{\prod}^\ast \in \mathcal{N}_{\mathcal{H}^q}$ are nuclear.
As a result of \eqref{equation of invertible operators in ARMA(p,q)} and \eqref{equation of invertible operators in ARMA(p,q), vector representation}, reasonable estimates for $\alpha_i$ and $\beta_j,$ respectively, are
\begin{gather}
  \hat{\alpha}_i  \coloneqq \hat{\psi}_i + \hat{B}_{[i]}\hat{\Psi}'''^{\,\top}_{[i]},\; 1\leq i\leq p, \qquad
  \hat{\beta}_j \coloneqq\begin{cases}
    \,\hat{B}^{(j)}_q, &\textnormal{if}\;\,1\leq j \leq q,\\
	\,0, &\textnormal{if}\;\,j > q.\\
  \end{cases}
\label{Estimators for alpha_i and beta_j in ARMA(p,q)}
\end{gather}
In these definitions, $\hat{B}_q$ stands for the estimator for $B_q$ defined by
\begin{gather*}
  \hat{B}_q \coloneqq - \hat{\Psi}'_{[q]}\,\boldsymbol{\hat{\prod}}^\dagger\coprod^{\hat{h}_{\!M}}_{\hat{h}_1},
\end{gather*}
where $\boldsymbol{\hat{\prod}}^\dagger\coloneqq(\boldsymbol{\hat{\prod}}\boldsymbol{\hat{\prod}}^\ast + \hspace{-0.05em}\gamma_N\hspace{0.05em}\mathbb{I})^{\hspace{-0.05em}-1},$ $(M_N) \subset\mathbb{N}$ and $ (\gamma_N)\subset(0,\infty)$ are sequences with $M_N\to \infty$ and $\gamma_N \rightarrow 0,$ and $\hat{h}_1, \dots, \hat{h}_{\!M}$  are the eigenfunctions associated to the eigenvalues $\hat{\zeta}_1\geq \cdots \geq\hat{\zeta}_M\geq 0$ of $\boldsymbol{\hat{\prod}}\boldsymbol{\hat{\prod}}^\ast.$ Further, $\hat{B}_{[i]},$ $B_{[i]},$ $ \hat{\Psi}'''_{[i]},$ $ \Psi'''_{[i]}\in \mathcal{S}_{\mathcal{H}^{\max(i,q)},\mathcal{H}}$ are for each $i,q$ defined by
\begin{gather}
  \hat{B}_{[i]} \coloneqq \begin{cases}
    \,\hat{B}_i, &\textnormal{if}\;\,1\leq i < q,\\
	\,\hat{B}_q, &\textnormal{if}\;\,i \geq q,\\
  \end{cases}\qquad\text{and}\qquad
  B_{[i]} \coloneqq
  \begin{cases}
    \,B_i, &\textnormal{if}\;\,1\leq i < q,\\
	\,B_q, &\textnormal{if}\;\,i \geq q,\\
  \end{cases}\label{Definition of modified vector of betas for estimation of alphas}\\
    \hat{\Psi}'''_{[i]} \coloneqq \begin{cases}
    \,\hat{\Psi}''_{[0]}, &\textnormal{if}\;\,1\leq i \leq q,\\
	\,\hat{\Psi}''_{[i-q]}, &\textnormal{if}\;\,i > q,\\
  \end{cases}\qquad\text{and}\qquad
  \Psi'''_{[q]}\coloneqq \begin{cases}
    \,\Psi''_{[0]}, &\textnormal{if}\;\,1\leq i \leq q,\\
	\,\Psi''_{[i-q]}, &\textnormal{if}\;\,i > q,\\
  \end{cases}\label{Definition of Psi'''}
\end{gather}
with $\hat{\Psi}''_{[i]}$ and $\hat{\Psi}'_{[i]}$ matching $\hat{\Psi}'_{[i]}$ in \eqref{definition of estimate for vector of psis for ARMA(p,q)} and $\Psi'_{[i]}$ in \eqref{definition of vector of psis for ARMA(p,q)} for all $i$, except for a reversed sign in the first component. To derive consistent estimates we make the following assumptions.
\begin{assumption}\label{As: operator for ARMA(p,q) estimation injective} The image of the operator-valued matrix $\boldsymbol{\prod}\colon\mathcal{H}^q\rightarrow \mathcal{H}^q$ lies dense.
\end{assumption}

 \begin{assumption}\label{As: Sobolev cond ARMA(p,q)} The operators $\alpha_1, \dots, \alpha_p, \beta_1, \dots, \beta_q\in\mathcal{S_H}$ are H-S, and for some $\gamma >0,$ it holds
  \begin{align*}
  \sum_{i=1}^\infty\sum_{j=1}^\infty\,\langle B_q, h_i\otimes d_j\rangle^2_\mathcal{S}\,(1 + i^{2\gamma}\! + \hspace{-0.05em}j^{2\gamma}) <
	\infty,
  \end{align*}
where $h_i$ and $d_j$ denote the eigenfunctions of $\boldsymbol{\prod}\boldsymbol{\prod}^\ast\colon \mathcal{H}^q \rightarrow \mathcal{H}^q$ and $\mathscr{C}_{\!\boldsymbol{X}}\colon \mathcal{H} \rightarrow \mathcal{H},$ respectively. 
\end{assumption}

\begin{assumption}\label{As: Dim and eigvals fARMA pq} There exists $d \in \mathbb{N}$ such that:
\begin{itemize}
    \item[\textnormal{(a)}] The dimensions of all eigenspaces of $\boldsymbol{\prod}\boldsymbol{\prod}^\ast$ are bounded above by $d;$
    \item[\textnormal{(b)}] For all $N$ sufficiently large, with probability one the eigenvalues of $\boldsymbol{\hat{\prod}}\boldsymbol{\hat{\prod}}^\ast$ satisfy $\hat{\zeta}_j\neq \hat{\zeta}_{j+1}$ for all $j\in\{1,...,M+d\}.$
\end{itemize}
\end{assumption}

\begin{theorem}\label{Theo: ARMA(p,q), infinite-dim} Suppose $(X_k)$ is the fARMA$(p,q)$ process in \eqref{functional ARMA definition}. Further, let the conditions of Theorem \ref{Theo: inv. op. - infinite-dim} and Assumptions \ref{As: operator for ARMA(p,q) estimation injective}–\ref{As: Dim and eigvals fARMA pq} hold, and let $M = M_{\!N} \to \infty$ with $M^{\gamma+1/2}Z^2_{\!M} = \mathrm{O}(K^\beta)$ and $\gamma_N = \mathrm O(\zeta_MM^{-\gamma})$, where $Z_{\!M} \coloneqq (\zeta_M - \zeta_{M_b})^{-1}$ as in \eqref{eq:Def Lambda_K}, with eigenvalues $\zeta_j$ of $\boldsymbol{\prod}\boldsymbol{\prod}^\ast$. Then,
\begin{align*}
  \max_{1\leq i \leq p, 1 \leq j \leq q}\!\big\{\|\hat{\alpha}_i - \alpha_i\|_\mathcal{S}, \|\hat{\beta}_j - \beta_j\|_\mathcal{S}\big\} &= \mathrm O_{\Prob}(M^{-\gamma}).
  \end{align*}
\end{theorem}

\section{Concluding remarks}\label{Section: Summary}

This article establishes consistent Yule--Walker-type estimators for the operators of invertible processes in separable Hilbert spaces under Sobolev conditions using Tikhonov regularization. We derive consistent estimators for functional linear, AR, MA, and ARMA processes of arbitrary order, leveraging asymptotic consistency results for growing finite-dimensional projections. The convergence rates primarily depend on operator smoothness in the inverted representation and error covariance regularity, with lesser sensitivity to operator norm residuals. Notably, high smoothness give near-parametric rate $N^{-1/2}$ even under rough error covariance operators.

\cite{GrangerAndersen1978} introduced a broader notion of invertibility for real-valued processes with i.i.d. errors, defining $(X_k)$ as invertible if $\Erw(\hat{\varepsilon}_k - \varepsilon_k)^2 \to 0$ as $k \to \infty$, where
\begin{align*}
  \varepsilon_k &= X_k -  f(X_{k-1}, \dots, X_{k-p}, \varepsilon_{k-1}, \dots, \varepsilon_{k-q}),\\
  \hat{\varepsilon}_k &= X_k -  f(X_{k-1}, \dots, X_{k-p}, \hat{\varepsilon}_{k-1}, \dots, \hat{\varepsilon}_{k-q}),
\end{align*}
with given initial values and i.i.d. centered innovations $\varepsilon_k$ with finite second moments. When $f$ is unknown, it is replaced by an estimate. \cite{Hallin1980, Hallin1981} extended this concept by replacing time limits with backward estimates, allowing for time dependence. This notion generalizes to separable Banach spaces $(\mathcal{B}, \|\cdot\|)$. A process $(X_k) \subset \mathcal{B}$ is \emph{mildly Granger--Andersen invertible} if
\begin{align*}
  \|\hat{\varepsilon}_{k_N}\! - \varepsilon_{k_N}\| = \mathrm o_{\Prob}(1),
\end{align*}
where $(\varepsilon_k)$ is i.i.d., $N$ is the sample size, and $k_N \to \infty$. This condition is weaker than standard invertibility \eqref{identity of invertible process}, requiring only convergence in probability without assuming finite second moments or constraining $(k_N)$. Standard invertibility implies mild Granger--Andersen invertibility, but not vice versa. As noted by \cite{GrangerAndersen1978}, certain functional time series, such as bilinear processes in separable Hilbert spaces \citep{hormann:kokoszka:2010}, satisfy mild but not standard invertibility. Others, like the non-linear MA$(1)$ process $X_k = \alpha_1(\varepsilon^2_{k-1}) + \varepsilon_k$ in $L^2[0,1]$, may be stationary yet fail to be invertible in both senses.

A natural direction for future research is extending estimation methods to mild Granger-Andersen invertibility under suitable assumptions. Further, developing a general estimation framework for Banach spaces, potentially even metric spaces, remains an open and promising challenge.

\appendix

\section{Appendix}\label{Section Proofs}

Various proofs utilize the \emph{operator-valued H\"older  inequality} for \emph{Schatten class operators}, along with several other standard results in operator theory \citep[][section 7.1]{Weidmann1980}. Let $\mathcal{H}_1$ and $\mathcal{H}_2$ be Hilbert spaces. Then, a compact linear operator $A\colon\mathcal{H}_1\rightarrow\mathcal{H}_2$ is a \emph{Schatten class operator of order} $1\leq p < \infty$ if its singular values $s_1(A)\geq s_2(A) \geq \dots \geq 0$ are $p$-fold summable. The set of such operators, denoted by $\mathcal{S}_{p,\mathcal{H}_1,\mathcal{H}_2}$ can be equipped with the norm $
\|A\|_p\coloneqq (\sum^\infty_{i=1}s^p_i(A))^{1/p},$ $A\in\mathcal{S}_{p,\mathcal{H}_1,\mathcal{H}_2}.$ The spaces of nuclear $\mathcal{N}_{\mathcal{H}_1,\mathcal{H}_2}$ and H-S operators $\mathcal{S}_{\mathcal{H}_1,\mathcal{H}_2}$ match $\mathcal{S}_{1,\mathcal{H}_1,\mathcal{H}_2}$ and $\mathcal{S}_{2,\mathcal{H}_1,\mathcal{H}_2},$ respectively. For the space of bounded linear operators $\mathcal{L}_{\mathcal{H}_1,\mathcal{H}_2},$ we write $\mathcal{S}_{\infty,\mathcal{H}_1,\mathcal{H}_2},$ and put $\|\cdot\|_\infty\coloneqq\|\cdot\|_\mathcal{L}.$ Note that for any $1\leq p < q \leq \infty$ it holds that $\mathcal{S}_{q,\mathcal{H}_1,\mathcal{H}_2} \subsetneq \mathcal{S}_{p,\mathcal{H}_1,\mathcal{H}_2},$ with $\|\cdot\|_p \leq \|\cdot\|_q.$

\begin{lemma}\label{Lemma: operator-valued hoelder inquality} Let $\mathcal{H}_1, \mathcal{H}_2, \mathcal{H}_3$ be Hilbert spaces, and let $p,q,r \in [1,\infty]$ with $\frac{1}{p} + \frac{1}{q} = \frac{1}{r},$ and $\frac{1}{\infty}\coloneqq 0.$ Then, for any $A\in\mathcal{S}_{p,\mathcal{H}_2,\mathcal{H}_3}$ and $B\in\mathcal{S}_{q,\mathcal{H}_1,\mathcal{H}_2},$ it holds that $AB\in\mathcal{S}_{r,\mathcal{H}_1,\mathcal{H}_3},$ with
$$
  \|AB\|_r \leq \|A\|_p\|B\|_q.
$$
\end{lemma}

\begin{proof}[\textbf{Proof}] See \cite{Weidmann1980}, Theorem 7.8.
\end{proof}

The following lemma follows from elementary norm inequalities and the details are omitted.

\begin{lemma}\label{Lemma: Inequalities norm of operator matrices} Let $\mathcal{H}$ be a Hilbert space, and let $B_{ij}\in\mathcal{L_H}$ and $S_{ij}\in\mathcal{S_H}$ be bounded and H-S operators, respectively, where $i=1, \dots, m\in\mathbb{N}, j=1, \dots, n\in\mathbb{N}.$ Then, the following holds.
\begin{itemize}
   \item[(a)] For the operator-valued matrix  $\boldsymbol{B} \coloneqq (B_{ij})_{1\leq i\leq m, 1\leq j\leq n}$ it holds that $\boldsymbol{B} \in \mathcal{L}_{\mathcal{H}^m\!, \mathcal{H}^n},$ with
  $$
  \|\boldsymbol{B}\|_{\mathcal{L}} \leq \sum^m_{i=1}\|(B_{i1}\;\cdots\;B_{in})\|_{\mathcal{L}}.
  $$
  \item[(b)] For the operator-valued matrix  $\boldsymbol{S} \coloneqq (S_{ij})_{1\leq i\leq m, 1\leq j\leq n}$ it holds that $\boldsymbol{S} \in \mathcal{S}_{\mathcal{H}^m\!, \mathcal{H}^n},$ with
  $$
    \|\boldsymbol{S}\|^2_{\mathcal{S}} = \sum^m_{i=1}\sum^n_{j=1}\|S_{ij}\|^2_{\mathcal{S}}.
  $$
\end{itemize}
\end{lemma}

In several places, we  make use of the following upper bounds. From \eqref{Nuclear norm cov op}, the definition of the elements of $\boldsymbol{X}^{[L]} = (X^{[L]}_k)_{k\in\mathbb{Z}}$ in \eqref{process X in Cartesian space}, and weak stationarity of $\boldsymbol{X}=(X_k)_{k\in\mathbb{Z}}$, it follows that
\begin{gather}\label{equality norm of cov-op of cartesion process}
  \|\mathscr{C}_{\!\boldsymbol{X}^{[L]}}\|_{\mathcal{N}}\! = \Erw\!\|X^{[L]}_0\|^2 = L\Erw\!\|X_0\|^2.
\end{gather}
\noindent Thus, due to \eqref{Nuclear norm cross-cov op} and the Cauchy--Schwarz inequality, $\|\mathscr{C}^h_{\!\boldsymbol{X}^{[L]}}\|_{\mathcal{N}}  \leq \|\mathscr{C}_{\!\boldsymbol{X}^{[L]}}\|_{\mathcal{N}}$ for all $h\in\mathbb{Z},$ and 
\begin{gather}
    \|\mathscr{C}^h_{\!\boldsymbol{X}^{[L]}, \boldsymbol{X}}\|_{\mathcal{N}} \leq \sqrt{L}\Erw\!\|X_0\|^2,\quad h\in\mathbb{Z}. \label{inequaility lag-h-cross-cov op}
\end{gather}

\begin{proof}[\textbf{Proof of Proposition \ref{Prop: Injectivity of C^L_X}}]

First, suppose $L=1$. The covariance operator $\mathscr{C}_{\!\boldsymbol{X}^{[1]}} = \mathscr{C}_{\!\boldsymbol{X}}$ of the invertible process in \eqref{identity of invertible process} can be expressed in terms of the lagged cross-covariance operators $\mathscr{C}^h_{\!\boldsymbol{X}} = \mathbb{E}(X_0 \!\otimes\! X_h)$ for $h \in \mathbb{Z}$, and the error covariance operator $\mathscr{C}_{\boldsymbol{\varepsilon}}$ as:
\begin{align}\label{C_X expansion of inv process X}
  \mathscr{C}_{\!\boldsymbol{X}} = \sum_{i=1}^{\infty} \sum_{j=1}^{\infty} \psi_j \mathscr{C}^{i-j}_{\!\boldsymbol{X}} \psi^*_i + \mathscr{C}_{\boldsymbol{\varepsilon}}.
\end{align}
The series converges in the nuclear norm topology, as the triangle inequality, operator-valued H\"older inequality, $\|\psi^\ast_i\|_\mathcal{L}=\|\psi_i\|_\mathcal{L},$ $\|\mathscr{C}^h_{\!\boldsymbol{X}}\|_{\mathcal{N}} \leq \|\mathscr{C}_{\!\boldsymbol{X}}\|_{\mathcal{N}}$ and Assumption \ref{As: Inv process} imply
\begin{align*}
  \Big\| \sum_{i=1}^{\infty} \sum_{j=1}^{\infty} \psi_j \mathscr{C}^{i-j}_{\!\boldsymbol{X}} \psi^\ast_i \Big\|_{\mathcal{N}} &\leq \|\mathscr{C}_{\!\boldsymbol{X}}\|_{\mathcal{N}} \Big( \sum_{i=1}^{\infty} \|\psi_i\|_{\mathcal{L}} \Big)^2<\infty.
\end{align*}
Since $\mathscr{C}_{\boldsymbol{\varepsilon}}$ is positive definite by Assumption \ref{As: Features of (eps_k)}, $\mathscr{C}_{\!\boldsymbol{X}}$ is injective if the double sum in \eqref{C_X expansion of inv process X} is positive semi-definite. In fact, for every $x \in \mathcal{H}$, 
\begin{align*}
  \Big\langle \sum_{i=1}^{\infty} \sum_{j=1}^{\infty} \psi_j \mathscr{C}^{i-j}_{\!\boldsymbol{X}} \psi^*_i (x), x \Big\rangle 
  &= \mathbb{E} \Big( \sum_{i=1}^{\infty}\, \langle X_i, \psi^*_i (x) \rangle \Big)^2 \geq 0.
\end{align*}

Next, let $L\in\mathbb{N}.$ Suppose by way of contradiction that there exists a $v = (v_1, \dots,v_L)^\top \in \mathcal{H}^L$ such that $v \ne 0,$ and
$$
0=\langle \mathscr{C}_{\!\boldsymbol{X}^{[L]}}(v), v \rangle =\! \mbox{ var}\Big(\sum_{i=0}^{{L-1}}\langle X_{-i}, v_i \rangle \Big).$$
If this is so there exists an $r \le L-1$ so that $v_r \ne 0$, and almost surely
\begin{align}\label{prop b proof}
\langle X_{-r}, v_r \rangle = -\sum_{i=r+1}^{L-1} \langle X_{-i}, v_i \rangle.
\end{align}
Notice the right-hand side of the above is measurable with respect to $\sigma( \varepsilon_i, \; i \le -r-1)$. Due to the invertible process representation \eqref{identity of invertible process}, the left-hand side of the above is equal to
$$
\sum^\infty_{j=1}\,\langle \psi_j(X_{-r-j}), v_r \rangle + \langle \varepsilon_{-r}, v_r \rangle.
$$
Therefore, multiplying each side of \eqref{prop b proof} with $\langle \varepsilon_{-r}, v_r \rangle$ and taking expectations gives $\langle \mathscr{C}_{\varepsilon}(v_r), v_r \rangle = 0$, with $v_r \ne 0$, which contradicts the injectivity of $\mathscr{C}_{\boldsymbol{\varepsilon}}.$
\end{proof}

To prove Theorem \ref{Theo: inv. op. - infinite-dim}, we first consider the idealized case where $\Psi_{\!L}$ is finite-dimensional (hence H-S) and diagonalizable with respect to the CONS $(c_i \otimes d_j) = (c_i(L) \otimes d_j)_{i,j}$ of $\mathcal{S}_{\mathcal{H}^L, \mathcal{H}}$, where $(\lambda_j, c_j) = (\lambda_j(L), c_j(L))_j$ and $(\mu_j, d_j)$ are the eigenpair sequences of $\mathscr{C}_{\!\boldsymbol{X}^{[L]}}$ and $\mathscr{C}_{\!\boldsymbol{X}}$, respectively.

\begin{assumption}\label{As: projection on finite-dim space for Psi_L} For   $K \in \mathbb{N}$ and  all $L \in \mathbb{N}$, there exist constants $p_{ijL}$ with $\sup_{i,j,L} p_{ijL}^2 < \infty$ such that
\begin{align*}
    \Psi_{\!L} = \Psi_{\!L}(K) &\coloneqq \sum_{i=1}^K\sum_{j=1}^K\,p_{ijL}(c_i\otimes d_j).
\end{align*}
\end{assumption}

\begin{lemma}\label{Lemma: Invertible operator estimation - finite-dim} Let Assumptions  \ref{As: Features of (eps_k)}--\ref{As: cons rate for lagged cov and cross-cov ops}, \ref{As: Dim eigenspaces unif bounded inv op}--\ref{As:No rep eigenvalues main Theo} hold. Further, suppose 
$$
  \Lambda_KL^{3/2} = \mathrm{o}(N^{1/2}), \quad \sum_{\ell>L}\|\psi_\ell\|_\mathcal{L} = \mathrm{O}(\lambda_K\Lambda^2_KL^2N^{-1/2}), \quad \theta_{\!N}=\mathrm O(\lambda_K\Lambda^2_KL^{5/2}N^{-1/2}).
$$
Then,
$$
    \|\hat{\Psi}_{\!L} - \Psi_{\!L}(K)\|_\mathcal{S} = \mathrm O_{\Prob}\big(\Lambda^2_KL^{5/2}N^{-1/2}\big).
$$
\end{lemma}

\begin{proof}[\textbf{Proof}] Note that $(\hat{\lambda}_j, \hat{c}_j) = (\hat{\lambda}_j(L), \hat{c}_j(L))_j$ is the eigenpair sequence of $\hat{\mathscr{C}}_{\!\boldsymbol{X}^{[L]}}$, and $\boldsymbol{X}^{[L]} = (X^{[L]}_k)_k \subset \mathcal{H}^L$ is the process in \eqref{process X in Cartesian space}. By Proposition \ref{Prop: Injectivity of C^L_X}, all eigenvalues $\lambda_j=\lambda_j(L)$ of $\mathscr{C}_{\!\boldsymbol{X}^{[L]}}$ are positive. Hereinafter, we write $\hat{\mathscr{D}}=\hat{\mathscr{D}}_{\!\boldsymbol{X}^{[L]}, \boldsymbol{X}},$ $\mathscr{D}=\mathscr{D}_{\!\boldsymbol{X}^{[L]}, \boldsymbol{X}},$ $\hat{\mathscr{C}}=\hat{\mathscr{C}}_{\!\boldsymbol{X}^{[L]}}$ and $\mathscr{C}=\mathscr{C}_{\!\boldsymbol{X}^{[L]}}.$ The definitions of $\hat{\Psi}_{\!L}$ and $\hat{\mathscr{C}}^\dagger = (\hat{\mathscr{C}} + \theta_{\!N} \mathbb{I})^{-1}$, the approximate Yule--Walker equation \eqref{approximate YW equation}, and $\mathscr{C}^\ddagger \coloneqq \mathscr{C} \mathscr{C}^\dagger$ yield
\begin{align}
  \hat{\Psi}_{\!L} - \Psi_{\!L}
  &= (\hat{\mathscr{D}} - \mathscr{D})\hat{\mathscr{C}}^\dagger\!\coprod^{\hat{c}_K}_{\hat{c}_1} + \,\mathscr{D}\Big(\hat{\mathscr{C}}^\dagger\!\coprod^{\hat{c}_K}_{\hat{c}_1} - \,\mathscr{C}^\dagger\!\coprod^{c_K}_{c_1}\Big) + \Big(\sum_{\ell>L}\psi_\ell\mathscr{C}^{1-\ell}_{\!\!\boldsymbol{X}^{[L]}, \boldsymbol{X}}\Big)\mathscr{C}^\dagger\!\coprod^{c_K}_{c_1} +  \,\Psi_{\!L}\Big(\mathscr{C}^\ddagger\!\coprod^{c_K}_{c_1} - \,\mathbb{I}\Big).\notag
\end{align}
The triangle inequality and the operator-valued H\"older inequality in Lemma \ref{Lemma: operator-valued hoelder inquality} give
\begin{align}
  \begin{split}\label{initial inequality operators of linear processes}
    \|\hat{\Psi}_{\!L} - \Psi_{\!L}\|_\mathcal{S} &\leq \| \hat{\mathscr{D}} - \mathscr{D}\|_\mathcal{S}\Big\|\hat{\mathscr{C}}^\dagger\!\coprod^{\hat{c}_K}_{\hat{c}_1}\Big\|_\mathcal{L} + \|\mathscr{D}\|_\mathcal{S}\Big\|\hat{\mathscr{C}}^\dagger\!\coprod^{\hat{c}_K}_{\hat{c}_1} - \,\mathscr{C}^\dagger\!\coprod^{c_K}_{c_1}\Big\|_\mathcal{L}\\
    &\qquad + \Big\|\sum_{\ell>L}\psi_\ell\mathscr{C}^{1-\ell}_{\!\!\boldsymbol{X}^{[L]}, \boldsymbol{X}}\Big\|_\mathcal{S}\Big\|\mathscr{C}^\dagger\!\coprod^{c_K}_{c_1}\Big\|_\mathcal{L} + \Big\|\Psi_{\!L}\Big(\mathscr{C}^\ddagger\!\coprod^{c_K}_{c_1}   - \;\mathbb{I}\Big)\Big\|_\mathcal{S}.  
  \end{split}
\end{align}
According to the definition of the operator norm, for all $K,L,N$, it holds
\begin{gather}
  \Big\|\hat{\mathscr{C}}^\dagger\!\coprod^{\hat{c}_K}_{\hat{c}_1}\Big\|_\mathcal{L} = \Big\|(\hat{\mathscr{C}} + \hspace{-0.05em}\theta_{\!N}\hspace{0.05em}\mathbb{I})^{\hspace{-0.05em}-1}\coprod^{\hat{c}_K}_{\hat{c}_1}\Big\|_\mathcal{L} = \sup_{j \leq K}(\hat{\lambda}_j + \theta_{\!N})^{-1} = (\hat{\lambda}_K + \theta_{\!N})^{-1},\label{Tikhonov est proj in inv}\\
    \Big\|\mathscr{C}^\dagger\!\coprod^{c_K}_{c_1}\Big\|_\mathcal{L} = \Big\|(\mathscr{C} + \hspace{-0.05em}\theta_{\!N}\hspace{0.05em}\mathbb{I})^{\hspace{-0.05em}-1}\coprod^{c_K}_{c_1}\Big\|_\mathcal{L} = \sup_{j \leq K}(\lambda_j + \theta_{\!N})^{-1} = (\lambda_K + \theta_{\!N})^{-1}.\label{Tikhonov pop proj in inv}
\end{gather}

With $K_b$ in the definition of $\Lambda_K=(\lambda_K-\lambda_{K_b})^{-1}$ in \eqref{eq:Def Lambda_K}, we have
\begin{align*}
    \Big\|\hat{\mathscr{C}}^\dagger\!\coprod^{\hat{c}_K}_{\hat{c}_1} - \;\mathscr{C}^\dagger\!\coprod^{c_K}_{c_1}\Big\|_\mathcal{L} 
    &\leq \Big\|\hat{\mathscr{C}}^\dagger\!\!\coprod^{\hat{c}_{K_b-1}}_{\hat{c}_1} - \;\mathscr{C}^\dagger\!\!\coprod^{c_{K_b-1}}_{c_1}\!\Big\|_\mathcal{L} + \Big\|\hat{\mathscr{C}}^\dagger\!\!\coprod^{\hat{c}_{K_b-1}}_{\hat{c}_{K+1}} - \;\mathscr{C}^\dagger\!\!\coprod^{c_{K_b-1}}_{c_{K+1}}\!\Big\|_\mathcal{L}\,.
\end{align*}
Further, by the definition of the given operators and the triangle inequality,
\begin{align}
  \begin{split}\label{Inequ for projected Tikh.}
     \Big\|\hat{\mathscr{C}}^\dagger\!\!\coprod^{\hat{c}_{K_b-1}}_{\hat{c}_1} - \;\mathscr{C}^\dagger\!\!\coprod^{c_{K_b-1}}_{c_1}\!\Big\|_\mathcal{L} &\leq \Big\|\sum^{K_b-1}_{i=1}\,(\hat{\lambda}_i+\theta_{\!N})^{-1}\Big[(\hat{c}_i\otimes \hat{c}_i) - (c_i\otimes c_i)\Big]\Big\|_\mathcal{L}\\
     &\qquad\; + \Big\|\sum^{K_b-1}_{i=1}\,\Big[(\hat{\lambda}_i+\theta_{\!N})^{-1} -(\lambda_i+\theta_{\!N})^{-1}\Big](c_i\otimes c_i)\Big\|_\mathcal{L}.  
  \end{split}
\end{align}
\noindent Applying Lemmas 3.1--3.2 from \cite{Reimherr2015} to \eqref{Inequ for projected Tikh.}, and using $\|\cdot\|_\mathcal{L} \leq \|\cdot\|_\mathcal{S}$ along with $(\sum_k a_k)^{1/2} \leq \sum_k a_k^{1/2}$ for non-negative $a_k$, with the reciprocal eigengaps $\hat{\Lambda}_K \coloneqq (\hat{\lambda}_K - \hat{\lambda}_{K+1})^{-1},$ which are well-defined by Assumption \ref{As:No rep eigenvalues main Theo}, it holds
\begin{align}\label{eq:est error Tikhonov proj, main theorem, main part}
    \Big\|\hat{\mathscr{C}}^\dagger\!\!\coprod^{\hat{c}_{K_b-1}}_{\hat{c}_1} - \;\mathscr{C}^\dagger\!\!\coprod^{c_{K_b-1}}_{c_1}\!\Big\|_\mathcal{L}  &\leq 2(\hat{\lambda}_{K_b-1}+\theta_{\!N})^{-1}\sqrt{K_b}\,\|\hat{\mathscr{C}} - \mathscr{C}\|_\mathcal{L}\big(\hat{\Lambda}_{K_b-1} + \Lambda_{K_b-1} + (\lambda_{K_b-1} + \theta_{\!N})^{-1} + 1 \big)\notag\\
    &\leq 2\hat{\Lambda}_{K_b-1}\sqrt{K+dL}\,\|\hat{\mathscr{C}} - \mathscr{C}\|_\mathcal{L}\big(\hat{\Lambda}_{K_b-1} + 2\Lambda_K + 1 \big),
\end{align}
where $K_b \leq K + dL$ is given by Assumption \ref{As: Dim eigenspaces unif bounded inv op}, and $\Lambda_K=\Lambda_{K_b-1}$ by definition. Similarly, 
\begin{align}\label{eq:est error Tikhonov proj, main theorem, remainder}
    \Big\|\hat{\mathscr{C}}^\dagger\!\!\coprod^{\hat{c}_{K_b-1}}_{\hat{c}_{K+1}} - \;\mathscr{C}^\dagger\!\!\coprod^{c_{K_b-1}}_{c_{K+1}}\!\Big\|_\mathcal{L} &\leq 2\hat{\Lambda}_{K_b-1}\sqrt{K_b-K}\,\|\hat{\mathscr{C}} - \mathscr{C}\|_\mathcal{L}\big(\hat{\Lambda}_{K_b-1} + 2\Lambda_{K_b-1} + 1 \big)\notag\\
    &\leq 2\hat{\Lambda}_{K_b-1}\sqrt{dL}\,\|\hat{\mathscr{C}} - \mathscr{C}\|_\mathcal{L}\big(\hat{\Lambda}_{K_b-1} + 2\Lambda_K + 1 \big).
\end{align}
\noindent In the following, we deduce $\hat{\Lambda}_{K_b-1} = \mathrm{O}_{\Prob}(\Lambda_K).$ Since
\begin{align*}
    \sup_{j \in \mathbb{N}} |\hat{\lambda}_j - \lambda_j| \leq \|\hat{\mathscr{C}}_{\!\boldsymbol{X}^{[L]}} - \mathscr{C}_{\!\boldsymbol{X}^{[L]}}\|_{\mathcal{L}} = \mathrm{O}_{\Prob}(L^{3/2}N^{-1/2})
\end{align*}
by Lemma 4.2 in \cite{Bosq2000} together with \eqref{up b cov}, and $\lambda_{K_b-1}=\lambda_K,$ it holds
\begin{align*}
  \hat{\lambda}_{K_b-1} - \hat{\lambda}_{K_b} &= \lambda_K - \lambda_{K_b} + (\hat{\lambda}_{K_b-1} - \lambda_{K_b-1}) + (\lambda_{K_b} - \hat{\lambda}_{K_b})\\
  &= \Lambda^{-1}_K + \mathrm{O}_{\Prob}(L^{3/2}N^{-1/2}).
\end{align*}
Therefore, due to the assumption $\Lambda_KL^{3/2} = \mathrm{o}(N^{1/2}),$ in fact 
\begin{align*}
 \hat{\Lambda}_K = \Lambda_K\big[1 + \mathrm{O}_{\Prob}(\Lambda_K L^{3/2}N^{-1/2})\big]^{-1} = \Lambda_K\big[1 + \mathrm{o}_{\Prob}(1)\big]^{-1} = \mathrm{O}_{\Prob}(\Lambda_K).
\end{align*}

Further, the triangle and the operator-valued H\"older inequality, and \eqref{inequaility lag-h-cross-cov op} yield for the nuclear norm $\|\cdot\|_\mathcal{N}\geq \|\cdot\|_\mathcal{S}$ of the remainder of the sum above for any $L,$
  \begin{align}\label{inequ sum operators}
    \Big\|\sum_{\ell>L}\psi_\ell\mathscr{C}^{1-\ell}_{\!\!\boldsymbol{X}^{[L]}, \boldsymbol{X}}\Big\|_\mathcal{N}\leq \sqrt{L}\Erw\!\|X_0\|^2\sum_{\ell>L}\|\psi_\ell\|_\mathcal{L},
  \end{align}
\noindent where the series on the right is finite by Assumption \ref{As: Inv process}. Finally, from the definition of the H-S norm, with $\Psi_{\!L}=\Psi_{\!L}(K)$ in Assumption \ref{As: projection on finite-dim space for Psi_L}, $ \mathscr{C}^\ddagger=\mathscr{C}(\mathscr{C} + \hspace{-0.05em}\theta_{\!N}\hspace{0.05em}\mathbb{I})^{\hspace{-0.05em}-1},$ and the Kronecker delta $\delta_{jk},$
\begin{align}
  \Big\|\Psi_{\!L}\Big(\mathscr{C}^\ddagger\!\coprod^{c_K}_{c_1}   - \,\mathbb{I}\Big)\Big\|^2_\mathcal{S} 
  &= \sum^\infty_{i=1}\,\big(\mathds 1_{[1,K]}(i)\lambda_i(\lambda_i + \theta_{\!N})^{-1}   - \,1\big)^2\,\Big\|\sum_{k=1}^K\sum_{\ell=1}^K\,p_{k\ell L}(c_k\otimes d_\ell)(c_i)\Big\|^2\notag\allowdisplaybreaks\\
  &= \sum^K_{i=1}\,\big(\lambda_i(\lambda_i + \theta_{\!N})^{-1}   - \,1\big)^2\,\Big\|\sum_{k=1}^K\sum_{\ell=1}^K\,p_{k\ell L}\delta_{ik}\,d_\ell\Big\|^2\notag\\
  &\leq \theta^2_{\!N}(\lambda_K+\theta_{\!N})^{-2}\|\Psi_{\!L}(K)\|^2_\mathcal{S}\;.\label{Psi_L(K) on difference of pseudo identity to identity}
\end{align}
\noindent Thereby, Parseval's identity yields
\begin{align}\label{squared norm of limit of projected operators of inv representation}
  \lim_{L\rightarrow\infty}\!\|\Psi_{\!L}(K)\|^2_\mathcal{S}\leq K^2\sup_{i,j,L}p^2_{ijL} < \infty.
\end{align}

Altogether, plugging \eqref{Tikhonov est proj in inv}--\eqref{Tikhonov pop proj in inv} and \eqref{eq:est error Tikhonov proj, main theorem, main part}--\eqref{Psi_L(K) on difference of pseudo identity to identity} in \eqref{initial inequality operators of linear processes}, and using $\|\mathscr{D}\|_\mathcal{S} \leq \sqrt{L}\Erw\!\|X_0\|^2$ by \eqref{inequaility lag-h-cross-cov op},  
$(\lambda_K + \theta_{\!N})^{-1} \leq \lambda^{-1}_K$ and $(\hat{\lambda}_K + \theta_{\!N})^{-1} \leq \hat{\Lambda}_{K_b-1},$ for each $K, L, N,$ we obtain
\begin{align*}
  \|\hat{\Psi}_{\!L} - \Psi_{\!L}(K)\|_\mathcal{S} &\leq \hat{\Lambda}_{K_b-1}\| \hat{\mathscr{D}} - \mathscr{D}\|_\mathcal{S}
  + 2\Erw\!\|X_0\|^2\hat{\Lambda}_{K_b-1}\sqrt{(K+dL)L}\,\|\hat{\mathscr{C}} - \mathscr{C}\|_\mathcal{L}\,(\hat{\Lambda}_{K_b-1} + 2\Lambda_K + 1)\notag\\
  &\quad\;\, + \lambda^{-1}_K\sqrt{L}\Erw\!\|X_0\|^2\sum_{\ell>L}\|\psi_\ell\|_\mathcal{L} + \theta_{\!N}\lambda^{-1}_K\|\Psi_{\!L}(K)\|_\mathcal{S}.
\end{align*}
Thus, due to $\hat{\Lambda}_{K_b-1} = \mathrm{O}_{\Prob}(\Lambda_K),$ Assumption \ref{As: cons rate for lagged cov and cross-cov ops}, our requirements and \eqref{squared norm of limit of projected operators of inv representation}, the claim follows.
\end{proof}

\begin{proof}[\textbf{Proof of Theorem \ref{Theo: inv. op. - infinite-dim}}] Follows from $
\|\hat{\Psi}_{\!L} -  \Psi_{\!L}\|_\mathcal{S} \leq \|\hat{\Psi}_{\!L} -  \Psi_{\!L}(K)\|_\mathcal{S} + \|\Psi_{\!L}(K) - \Psi_{\!L}\|_\mathcal{S},
$ the arguments to prove Lemma \ref{Lemma: Invertible operator estimation - finite-dim}, $\|\Psi_{\!L}(K) - \Psi_{\!L}\|_\mathcal{S} = \mathrm{O}( K^{-\beta})$ and $\|\Psi_{\!L}\|^2_\mathcal{S} = \sum_{i=1}^\infty \|\psi_j\|^2_\mathcal{S} < \infty$ (Remark \ref{Remark on operators in inverted representation} (c)) by Assumption \ref{As: Sob cond inv ops}, along with the requirements of our sequences.
\end{proof}

\begin{proof}[\textbf{Proof of Proposition \ref{Prop: Unique operators in LP and inv process}}] (a)\, According to Fubini's theorem, and as linear operators commute with the expected value, we have for any $i,j,k,$
\begin{align}\label{eqln}
\Erw  \langle c_k, \varepsilon_{i-j} \rangle \sum_{\ell=0}^{\infty} A_\ell(\varepsilon_{i-\ell}) &= \sum_{\ell=0}^{\infty}A_\ell( \Erw \langle c_k, \varepsilon_{i-j} \rangle \varepsilon_{i-\ell}).
\end{align}
We note that $\Erw  \langle c_k, \varepsilon_{i-j} \rangle \varepsilon_{i-\ell}=0$ when $j\ne \ell$, and when $j=\ell$,
$$ \Erw  \langle c_k, \varepsilon_{i-j} \rangle \varepsilon_{i-j}= \sum_{r=0}^{\infty}\Erw \langle c_k, \varepsilon_{i-j} \rangle \langle c_r, \varepsilon_{i-j} \rangle c_r = \lambda_k c_k$$
after \cite{Bosq2000}, equation (1.36). Combining with  \eqref{eqln}, we have that
$$ 
  \Erw  \langle c_k, \varepsilon_{i-j} \rangle \sum_{\ell=0}^{\infty} A_\ell(\varepsilon_{i-\ell}) = \lambda_k A_j(c_k),\quad \mbox{ and } \quad
  \Erw  \langle c_k, \varepsilon_{i-j} \rangle \sum_{\ell=0}^{\infty} B_\ell(\varepsilon_{i-\ell}) = \lambda_k B_j(c_k).
$$
Hence, due to $\lambda_k A_j(c_k) = \lambda_k B_j(c_k)$ for all $j,k\in \mathbb{N},$ and as $\mathscr{C}_{\boldsymbol{\varepsilon}}$ is injective, $\lambda_k >0$, thus $A_j(c_k) =   B_j(c_k)$ for all $j,k\in \mathbb{N}.$ Consequently,  as $(c_j)$ forms a CONS of $\mathcal{H}$, for all $j$ holds indeed $A_j = B_j.$

\vspace{1ex}

\noindent (b)\, Note that if $(X_k)$ is a centered linear process, then $Z_i = \sum_{\ell=0}^\infty A_\ell(X_{i-\ell})$ is also one, since
$$
Z_i = \sum^\infty_{\ell=0}A_\ell(X_{i-\ell}) = \sum^\infty_{\ell=0}D_\ell(\varepsilon_{i-\ell}),
$$
where $D_\ell = \sum_{j=0}^\ell A_j\phi_{\ell - j}. $ By hypothesis we also have that
$$
Z_i = \sum^\infty_{\ell=0}B_\ell(X_{i-\ell}) = \sum^\infty_{\ell=0}F_\ell(\varepsilon_{i-\ell}),
$$
where $F_\ell = \sum_{j=0}^\ell B_j\phi_{\ell - j}. $ According to part (a), $F_\ell = D_\ell$ for all $\ell \in \mathbb{N}$. Since $\phi_0 = \mathbb{I}$  we have then that $A_0 = D_0  =  F_0 = B_0$. Now the result follows by induction.
\end{proof}

We now turn towards proving Proposition \ref{Proposition linear operator estimation - infinite-dim}. Under Assumption \ref{As: projection on finite-dim space for Psi_L} it follows for any $\ell$ and  $K$ and  $L$ sufficiently large that
\begin{align}\label{finite-dim representation of psi_l}
  \psi_\ell=\psi_\ell(K)\coloneqq \big[\Psi_{\!L}(K)\big]^{(\ell)} = \sum_{i=1}^K\sum_{j=1}^K\,p_{ijL}(c_i\otimes d_j)^{(\ell)}\, = \sum_{i=1}^K\sum_{j=1}^K\,p_{ijL}(c^{(\ell)}_i\!\otimes d_j).
\end{align}
The finite-dimensional representations of $\phi_i$ need not be straightforward. For instance, for $\phi_2$, using $\phi_1=\psi_1$ and $\phi_2 = \psi_1\phi_1$ under Assumption \ref{As: projection on finite-dim space for Psi_L} for sufficiently large $L$, we obtain
\begin{align}
\phi_2(K) = \psi^2_1(K) &= \sum_{i=1}^K\sum_{i'=1}^K\sum_{j=1}^K\sum_{j'=1}^K\,p_{ijL}\,p_{i'j'L}\underbrace{\big(c^{(\ell)}_i\!\otimes d_j\big)\big(c^{(\ell)}_{i'}\!\otimes d_{j'}\big)}_{=\;\langle c^{(\ell)}_i\!,\,d_{j'}\rangle\,c^{(\ell)}_{i'}\otimes\, d_j}.\notag
\end{align}
Nevertheless, we can use the recursively defined finite-dimensional approximations
\begin{align}\label{implicit finite-dim representation of phi_i}
\phi_i(K) \coloneqq  \sum^i_{j=1}\,\psi_j(K)\phi_{i-j}(K),\quad i\in\mathbb{N},\qquad\text{with}\;\, \phi_0(K)\coloneqq \mathbb{I}.
\end{align}

\begin{proposition}\label{Proposition linear operator estimation - finite-dim} Let the assumptions of Lemma \ref{Lemma: Invertible operator estimation - finite-dim} hold. Then, 
  \begin{align*}
    \|\hat{\phi}_i - \phi_i(K)\|_\mathcal{S}  = \mathrm{O}_{\Prob}\big(\|\hat{\Psi}_{\!L} - \Psi_{\!L}(K)\|_\mathcal{S}\big) = \mathrm O_{\Prob}\big(\Lambda^2_KL^{5/2}N^{-1/2}\big),\quad i,K\in\mathbb{N}.
  \end{align*}
\end{proposition}

\begin{proof}[\textbf{Proof}] Due to the definition of $\hat{\phi_i}$ in \eqref{Definition of estimators for phi_i} and $\phi_i(K)$ in \eqref{implicit finite-dim representation of phi_i}, and since for each $L\geq j$ and $K$, $\|\hat{\psi_j} - \psi_j(K)\|_\mathcal{S} \leq \|\hat{\Psi}_{\!L} - \Psi_{\!L}(K)\|_\mathcal{S},$  elementary conversions give for each $i\leq L$ and $K\colon$
\begin{align}
\|\hat{\phi}_i - \phi_i(K)\|_\mathcal{S} 
  &\leq \sum^{i-1}_{j=1}\,\|\hat{\psi}_j\|_\mathcal{L}\|\hat{\phi}_{i-j} - \phi_{i-j}(K)\|_\mathcal{S} + \|\hat{\Psi}_{\!L} - \Psi_{\!L}(K)\|_\mathcal{S}\sum^{i-1}_{j=0}\,\|\phi_j\|_\mathcal{L},\notag
\end{align}
where we used $ \hat{\phi}_0=\phi_0(K)=\mathbb{I}$ and $\|\phi_j(K)\|_\mathcal{L} \leq \|\phi_j\|_\mathcal{L}$. Since $\|\hat{\phi}_1 - \phi_1(K)\|_\mathcal{S} = \|\hat{\psi}_1 - \psi_1(K)\|_\mathcal{S}$ is bounded above by $\|\hat{\Psi}_{\!L} - \Psi_{\!L}(K)\|_\mathcal{S}$, whose asymptotic upper bound follows from Lemma \ref{Lemma: Invertible operator estimation - finite-dim}, iterating the above $i-1$ times yields an asymptotic bound for $\|\hat{\phi}_i - \phi_i(K)\|_\mathcal{S}$. Namely, with
\begin{align}\label{Theta i}
    \Theta(i) \coloneqq 1 + \sum^{i-1}_{k=1}\sum^{i-j_0-1}_{j_1=1}\cdots\sum^{i-j_{k-1}-1}_{j_k=1}\!\Big(\prod^{k-1}_{\ell=1}\,\|\psi_{j_k}\|_\mathcal{L}\!\Big)\Big[\,\|\phi_{j_k}\|_\mathcal{L} + \|\psi_{j_k}\|_\mathcal{L}\Big], \quad i\in\mathbb{N},
\end{align}
and $\|\phi_0\|_\mathcal{L}=1,$ $\|\hat{\psi}_j\|_\mathcal{L} = \|\psi_j(K)\|_\mathcal{L} + \mathrm o_{\Prob}(1)$ for all $j,K$ by Lemma \ref{Lemma: Invertible operator estimation - finite-dim} and $\|\psi_j(K)\|_\mathcal{L} \leq \|\psi_j\|_\mathcal{L},$ 
\begin{align}
  \|\hat{\phi}_i - \phi_i(K)\|_\mathcal{S} &\leq \|\hat{\Psi}_{\!L} - \Psi_{\!L}(K)\|_\mathcal{S}\bigg[\sum^{i'}_{j_1=0}\|\phi_{j_1}\|_\mathcal{L} + \sum^{i'}_{j_1=1}\sum^{i'-j_1}_{j_2=0}\|\hat{\psi}_{j_1}\|_\mathcal{L}\|\phi_{j_2}\|_\mathcal{L}\bigg.\notag\\
   &\bigg.\qquad\qquad\qquad\qquad\quad\;\; + \cdots + \sum^{i'}_{j_1=1}\sum^{i'-j_1}_{j_2=1}\cdots\sum^{i'-j_{i-3}}_{j_{i-2}=1}\sum^{i'-j_{i-2}}_{j_{i-1}=0}\!\Big(\prod^{i'-1}_{k=1}\|\hat{\psi}_{j_k}\|_\mathcal{L}\Big)\|\phi_{j_{i-1}}\|_\mathcal{L}\bigg]\notag\allowdisplaybreaks\\ 
   &\qquad\qquad + \sum^{i'}_{j_1=1}\sum^{i'-j_1}_{j_2=1}\cdots\sum^{i'-j_{i-2}}_{j_{i-1}=1}\!\Big(\prod^{i'}_{k=1}\|\hat{\psi}_{j_k}\|_\mathcal{L}\Big)\underbrace{\;\big\|\hat{\phi}_{i-\sum^{i-1}_{k=1}j_k}\! - \phi_{i-\sum^{i-1}_{k=1}j_k}\big\|_\mathcal{S}\;}_{=\;\|\hat{\psi}_1 - \psi_1(K)\|_\mathcal{S}} \notag\allowdisplaybreaks\\
 &\leq \|\hat{\Psi}_{\!L} - \Psi_{\!L}(K)\|_\mathcal{S}\big[\Theta(i) + \mathrm o_{\Prob}(1)\big],\label{upper bound of hat phi i - phi i K}
\end{align}
with $i' \coloneqq i-1$, where factors in multi-sums with $\sum^0_{j_k=1}$ for some $k$ were set to 1.
\end{proof}

\begin{proof}[\textbf{Proof of Proposition \ref{Proposition linear operator estimation - infinite-dim}}] 
Similar to the arguments in the proof of Proposition \ref{Proposition linear operator estimation - finite-dim}, and since $\|\psi_j(K) - \psi_j\|_\mathcal{S} \leq \|\Psi_{\!L}(K) - \Psi_{\!L}\|_\mathcal{S}$ for all $j, K,$ for $i\leq L,$ it holds that
\begin{align*}
  \|\phi_i(K) - \phi_i\|_\mathcal{S}   &\leq \sum^{i-1}_{j=1}\,\|\psi_j\|_\mathcal{L} \|\phi_{i-j}(K) - \phi_{i-j}\|_\mathcal{S} + \|\Psi_{\!L}(K) - \Psi_{\!L}\|_\mathcal{S} \sum^{i-1}_{j=0}\,\|\phi_j\|_\mathcal{L}.
\end{align*}
Performing $i-1$ iterations of this inequality, following the rationale to prove Proposition \ref{Proposition linear operator estimation - finite-dim}, and using $\|\phi_1(K) - \phi_1\|_\mathcal{S} = \|\psi_1(K) - \psi_1\|_\mathcal{S},$ gives $\|\phi_i(K) - \phi_i\|_\mathcal{S} \leq \|\Psi_{\!L}(K) - \Psi_{\!L}\|_\mathcal{S} \Theta(i),$ with $\Theta(i)$ in \eqref{Theta i}. Further, \eqref{upper bound of hat phi i - phi i K}, $\|\hat{\Psi}_{\!L} - \Psi_{\!L}(K)\|_\mathcal{S} = \mathrm{O}_{\Prob}(K^{-\beta})$ and $\|\Psi_{\!L}(K) - \Psi_{\!L}\|_\mathcal{S} = \mathrm{O}(K^{-\beta})$ due to the requirements of  Theorem \ref{Theo: inv. op. - infinite-dim}, and $\|\hat{\phi}_i - \phi_i\|_\mathcal{S} \leq 
  \|\hat{\phi}_i - \phi_i(K)\|_\mathcal{S} + \|\phi_i(K) - \phi_i\|_\mathcal{S}$ prove our claim. 
\end{proof}

\begin{proof}[\textbf{Proof of Corollary \ref{Corollary AR(p) - finite-dim}}] Follows directly from Theorem \ref{Theo: inv. op. - infinite-dim} with $L=p.$ 
\end{proof}

\begin{proof}[\textbf{Proof of Corollary \ref{Corollary from linear operators}}] Follows from Proposition \ref{Proposition linear operator estimation - infinite-dim} by putting $\beta_j = \phi_j$ and $\hat{\beta}_j = \hat{\phi}_j.$
\end{proof}

To derive the asymptotic result for our fARMA$(1,1)$ operators, we derive a preliminary result based on the following assumption. 

\begin{assumption}\label{As: projection on finite-dim space for alpha_1 and beta_1} For fixed $M\in\mathbb{N},$ there are constants $b_{ij}\in\mathbb{R}$ so that 
  \begin{align*}
    \beta_1 = \beta_1(M) &\coloneqq \sum_{i=1}^M\sum_{j=1}^M\,b_{ij}(f_i\otimes f_j).
  \end{align*}
\end{assumption}

\noindent A finite-dimensional representation of $\beta_1$ does not imply one for $\alpha_1$ under Assumption \ref{As: psi_1psiadj_1 injective}. As \eqref{equation of invertible operators in ARMA(1,1)} gives $\alpha_1 = \alpha_1(M) \coloneqq \psi_1 + \beta_1(M)$, assuming $\psi_1$ is finite-dimensional contradicts Assumption \ref{As: psi_1psiadj_1 injective}.

\begin{lemma}\label{Lemma for ARMA(1,1) operators, finite-dim} 
Suppose $(X_k)$ is the fARMA$(1,1)$ process in \eqref{fARMA(1,1) definition}. Further,  let the conditions of Theorem \ref{Theo: inv. op. - infinite-dim}, Assumptions \ref{As: psi_1psiadj_1 injective}, \ref{As: Dim and eigvals fARMA 11}, \ref{As: projection on finite-dim space for alpha_1 and beta_1}, and $\gamma_N = \mathrm{O}(K^{-\beta})$ hold. Then, 
\begin{align*}
    \max\!\big\{\|\hat{\alpha}_1 - \alpha_1(M)\|_\mathcal{S}, \|\hat{\beta}_1 - \beta_1(M)\|_\mathcal{S}\big\} &= 
    \mathrm O_{\Prob}\big(\|\hat{\Psi}_{\!L} - \Psi_{\!L}\|_\mathcal{S}\big) = 
    \mathrm O_{\Prob}(K^{-\beta}).
\end{align*}
\end{lemma}

\begin{proof}[\textbf{Proof}] The conversions here follow those in Lemma \ref{Lemma: Invertible operator estimation - finite-dim}. Recall $(\gamma_N) \subset (0,\infty)$ with $\gamma_N \to 0$ and $\beta_1 = \beta_1(M)$ by  Assumption \ref{As: projection on finite-dim space for alpha_1 and beta_1}. Let $(\rho_j, f_j)$ and $(\hat{\rho}_j, \hat{f}_j)$ be the eigenpair sequences of $\psi_1\psi_1^\ast$ and $\hat{\psi}_1\hat{\psi}_1^\ast$, respectively, where $\psi_i$ are the H-S operators from \eqref{identity of invertible process}, and $\hat{\psi}_1 = \hat{\Psi}^{(i)}_{\!L} $ its estimates from \eqref{Operator estimate for big Psi_L}. Since $\psi_1\psi_1^\ast$ is injective (Assumption \ref{As: psi_1psiadj_1 injective}), all $\rho_j$ are positive. Moreover,
\begin{align*}
  \|\hat{\psi}_1\hat{\psi}^\ast_1 - \psi_1\psi^\ast_1\|_\mathcal{L} \leq \|\hat{\psi}_1 - \psi_1\|_\mathcal{L}\big(\|\psi_1\|_\mathcal{L} + \|\hat{\psi}_1\|_\mathcal{L}\big).
\end{align*}

As $\hat{\beta}_1 \coloneqq - \hat{\psi}_2\hat{\psi}^\dagger_1\coprod^{\hat{f}_M}_{\hat{f}_1},$ where $ \hat{\psi}^\dagger_1 =  \hat{\psi}^\ast_1(\hat{\psi}_1\hat{\psi}^\ast_1 + \gamma_N\mathbb{I})^{-1}\!,$ with  $\psi^\dagger_1 \coloneqq  \psi^\ast_1(\psi_1\psi^\ast_1 + \gamma_N\mathbb{I})^{-1}\!, \psi^\ddagger_1 \coloneqq \psi^\ast_1\psi^\dagger_1,$ and  $\psi_2 = -\beta_1\psi_1$ from \eqref{equation of invertible operators in ARMA(1,1)}, due to triangle and  operator-valued H\"older's inequality, it holds
\begin{align*}
  \|\hat{\beta}_1  - \beta_1\|_\mathcal{S} \leq \|\hat{\psi}_2 - \psi_2\|_\mathcal{S}\Big\|\hat{\psi}^\dagger_1\coprod^{\hat{f}_M}_{\hat{f}_1}\Big\|_\mathcal{L} + \|\psi_2\|_\mathcal{S}\Big\|\hat{\psi}^\dagger_1\coprod^{\hat{f}_M}_{\hat{f}_1} - \psi^\dagger_1\coprod^{f_M}_{f_1}\Big\|_\mathcal{L} + \Big\|\beta_1\Big(\psi^\ddagger_1\coprod^{f_M}_{f_1} - \mathbb{I}\Big)\Big\|_\mathcal{S}\;.
\end{align*}
\noindent Due to the properties of our operators and their norms, it follows that
\begin{align*}
  \Big\|\hat{\psi}^\dagger_1\coprod^{\hat{f}_M}_{\hat{f}_1}\Big\|^2_\mathcal{L} &= \Big\|\hat{\psi}_1\hat{\psi}^\ast_1\big(\hat{\psi}_1\hat{\psi}^\ast_1 + \gamma_N\mathbb{I}\big)^{-2}\coprod^{\hat{f}_M}_{\hat{f}_1}\Big\|_\mathcal{L}\leq (\hat{\rho}_M + \gamma_N)^{-1}.
\end{align*}
\noindent Further, from the definition of $\hat{\psi}^\dagger_1$ and $\psi^\dagger_1,$ since $\|A^\ast\|_\mathcal{L} = \|A\|_\mathcal{L}$ for bounded linear operators $A,$ elementary transformations, and similar  arguments as in the proof of Lemma  \ref{Lemma: Invertible operator estimation - finite-dim}, it follows
\begin{align*}
  &\Big\|\hat{\psi}^\dagger_1\coprod^{\hat{f}_M}_{\hat{f}_1} - \,\psi^\dagger_1\coprod^{f_M}_{f_1}\Big\|_\mathcal{L}\notag\\
  &\;\leq \|\hat{\psi}_1\|_\mathcal{L}\Big\|(\hat{\psi}_1\hat{\psi}^\ast_1 + \gamma_N\mathbb{I})^{-1}\coprod^{\hat{f}_M}_{\hat{f}_1} - \,(\psi_1\psi^\ast_1 + \gamma_N\mathbb{I})^{-1}\coprod^{f_M}_{f_1}\Big\|_\mathcal{L} + \|\hat{\psi}_1 - \psi_1\|_\mathcal{L}\Big\|(\psi_1\psi^\ast_1 + \gamma_N\mathbb{I})^{-1}\coprod^{f_M}_{f_1}\Big\|_\mathcal{L}
  \notag\\
  &\;\leq 2\hat{P}_{M_b-1}\sqrt{M+d}\,\|\hat{\psi}_1\hat{\psi}^\ast_1 - \psi_1\psi^\ast_1\|_\mathcal{L}\big(\hat{P}_{M_b-1} + 2P_M + 1 \big) + \rho^{-1}_{\!M} \|\hat{\psi}_1 - \psi_1\|_\mathcal{L}, 
\end{align*}
with $\hat{P}_{M_b-1}\coloneqq (\hat{\rho}_{M_b-1} - \hat{\rho}_{M_b})^{-1},$ and  $P_{M}\coloneqq (\rho_M - \rho_{M_b})^{-1}$ as in \eqref{eq:Def Lambda_K}, where former is well-defined by Assumption \ref{As: Dim and eigvals fARMA 11}. We also used that $M_b \leq M + d$ for some $d\in\mathbb{N}$ by Assumption \ref{As: Dim and eigvals fARMA 11}. Moreover, as in the proof of Lemma \ref{Lemma: Invertible operator estimation - finite-dim}, with $\beta_1=\beta_1(M)$ in Assumption \ref{As: projection on finite-dim space for alpha_1 and beta_1}, it holds
\begin{align*}
  \Big\|\beta_1\Big(\psi^\ddagger_1\coprod^{f_M}_{f_1} - \mathbb{I}\Big)\Big\|_\mathcal{S} 
  &\leq \gamma_N\rho^{-1}_{\!M}\|\beta_1(M)\|_\mathcal{S}.
\end{align*}

Altogether, with $\|\cdot\|_\mathcal{L} \leq \|\cdot\|_\mathcal{S},$ $\|\hat{\psi}_i - \psi_i\|_\mathcal{S} \leq \|\hat{\Psi}_{\!L} - \Psi_{\!L}\|_\mathcal{S}$ for all $i$, and $(\hat{\rho}_M + \gamma_N)^{-1} \leq \hat{P}_{M_b-1},$ $\hat{P}_{M_b-1} = \mathrm O_{\Prob}(P_M)$ (see the proof of Lemma \ref{Lemma: Invertible operator estimation - finite-dim}), $P^{1/2}_M = \mathrm{O}(P_M),$ and $\rho^{-1}_M \leq P_M,$ it holds 
\begin{align}
    \|\hat{\beta}_1  - \beta_1(M)\|_\mathcal{S} 
    &\leq \mathrm{O}_{\Prob} \big(M^{1/2} P^2_M \|\hat{\Psi}_{\!L} - \Psi_{\!L}\|_\mathcal{S} \big) + \mathrm{O}(\gamma_N \rho^{-1}_{\!M})\label{upper bound beta_1 M fixed}.
\end{align}
Thus, for fixed $M,$  $\|\hat{\beta}_1  - \beta_1(M)\|_\mathcal{S} = \mathrm O_{\Prob}(K^{-\beta})$ by Theorem \ref{Theo: inv. op. - infinite-dim} and $\gamma_N = \mathrm{O}(K^{-\beta}).$ Consequently, due to $\hat{\alpha}_1 = \hat{\psi}_1 - \hat{\beta}_1$ and $\alpha_1(M) = \psi_1 - \beta_1(M),$ for fixed $M,$ we also obtain $\|\hat{\alpha}_1  - \alpha_1(M)\|_\mathcal{S}\!\leq\!\|\hat{\Psi}_{\!L} - \Psi_{\!L}\|_\mathcal{S} + \|\hat{\beta}_1 - \beta_1(M)\|_\mathcal{S} = \mathrm O_{\Prob}(K^{\-\beta}).$ 
\end{proof}

\begin{proof}[\textbf{Proof of Theorem \ref{Theorem for ARMA(1,1) operators, infinite-dim}}] Due to \eqref{upper bound beta_1 M fixed}, $M^{\gamma+1/2}P^2_M = \mathrm{O}(K^\beta)$, $\gamma_N = \mathrm{O}(\rho_M M^{-\gamma})$ and Assumption \ref{As: Sobolev cond ARMA(1,1) - alpha_1beta_1}, we have $\|\hat{\beta}_1 - \beta_1\|_\mathcal{S} = \mathrm{O}_{\Prob}(M^{-\gamma})$ (see the proof of Theorem \ref{Theo: inv. op. - infinite-dim}). In consequence, due to $\alpha_1 = \psi_1 + \beta_1,$ Theorem \ref{Theo: inv. op. - infinite-dim}, and $K^{-\beta} = \mathrm{O}(M^{-\gamma}),$ it also holds $\|\hat{\alpha}_1 - \alpha_1\|_\mathcal{S} = \mathrm{O}_{\Prob}(M^{-\gamma})$.
\end{proof}

\begin{assumption}\label{As: projection on finite-dim space for B_q in ARMA(p,q)} For fixed $M\in\mathbb{N},$ there are constants $b_{ijq}\in\mathbb{R}$ so that 
  \begin{align*}
    B_q = B_q(M) &\coloneqq \sum_{i=1}^M\sum_{j=1}^M\,b_{ijq}(h_i\otimes d_j).
  \end{align*}
\end{assumption}
\noindent From \eqref{equation of invertible operators in ARMA(p,q)}, which is equivalent to $\alpha_i = \psi_i + B_{[i]}\Psi'''^{\,\top}_{[i]}$ for $i=1, \dots, p,$ as in Section \ref{Subsec: fARMA}, it follows
\begin{align*}
 \alpha_i = \alpha_i(M) \coloneqq \psi_i + B_{[i]}(M)\Psi''^{\,\top}_{[i]},\quad 1\leq i \leq p.
\end{align*}

\begin{lemma}\label{Theorem for ARMA(p,q) operators, finite-dim}
Let $(X_k)$ be the fARMA$(p,q)$ process in \eqref{functional ARMA definition}. Moreover, let the conditions of Theorem \ref{Theo: inv. op. - infinite-dim}, Assumptions \ref{As: operator for ARMA(p,q) estimation injective}, \ref{As: Dim and eigvals fARMA pq}, \ref{As: projection on finite-dim space for B_q in ARMA(p,q)} and $\gamma_N = \mathrm{O}(K^{-\beta})$ hold.  Then, 
$$
  \max_{1\leq i \leq p, 1\leq j \leq q}\!\big\{\|\hat{\alpha}_i - \alpha_i(M)\|_\mathcal{S},  \|\hat{\beta}_j - \beta_j(M)\|_\mathcal{S}\big\} = \mathrm O_{\Prob}\big(\|\hat{\Psi}_{\!L} - \Psi_{\!L}\|_\mathcal{S}\big) = \mathrm O_{\Prob}(K^{-\beta}).
$$
\end{lemma}

\begin{proof}[\textbf{Proof}] The proof follows that of Lemma \ref{Lemma for ARMA(1,1) operators, finite-dim}. Let $(\zeta_j, h_j)$ and $(\hat{\zeta}_j, \hat{h}_j)$ be the eigenpair sequences of $\boldsymbol{\prod}\boldsymbol{\prod}^\ast, \boldsymbol{\hat{\prod}}\boldsymbol{\hat{\prod}}^\ast \in \mathcal{N}_{\mathcal{H}^q},$ with $\boldsymbol{\prod}, \boldsymbol{\hat{\prod}} \in \mathcal{S}_{\mathcal{H}^q}$ in \eqref{definition of matrices of psis for ARMA(p,q)} and \eqref{definition of estimate for matrices of psis for ARMA(p,q)}, respectively. Under Assumption \ref{As: Sob cond inv ops} and Lemma \ref{Lemma: Inequalities norm of operator matrices} (b), we have $\boldsymbol{\prod}\boldsymbol{\prod}^\ast \in \mathcal{N}_{\mathcal{H}^q},$ and Assumption \ref{As: operator for ARMA(p,q) estimation injective} ensures all eigenvalues $\zeta_j$ are positive. Finally, basic conversions and the definitions of our operator-valued vectors and matrices yield
\begin{align*}
  \Big\|\boldsymbol{\hat{\prod}} - \boldsymbol{\prod}\Big\|_\mathcal{L} &\leq \sum^{q-1}_{i=0}\,\|\hat{\Psi}'_{[i]} - \Psi'_{[i]}\|_\mathcal{L}\,.
\end{align*}
Consequently, as $\|\hat{\Psi}'_{[i]} - \Psi'_{[i]}\|_\mathcal{S}\leq \|\hat{\Psi}_{\!L} - \Psi_{\!L}\|_\mathcal{S}$ for $L\geq p+q+i-1,$  we obtain for $L\geq p+2q-2,$
\begin{align}
  \Big\|\boldsymbol{\hat{\prod}}\boldsymbol{\hat{\prod}}^\ast\! - \boldsymbol{\prod}\boldsymbol{\prod}^\ast\Big\|_\mathcal{L} 
  &\leq q\|\hat{\Psi}_{\!L} - \Psi_{\!L}\|_\mathcal{S}\Big(\Big\|\boldsymbol{\prod}\Big\|_\mathcal{L} + \Big\|\boldsymbol{\hat{\prod}}\Big\|_\mathcal{L}\Big).\label{inequality of difference of estimate to in ARMA(p,q)}
\end{align}

With $\hat{B}_q = - \hat{\Psi}'_{[q]}\boldsymbol{\hat{\prod}}^\dagger\!\coprod^{\hat{h}_{\!M}}_{\hat{h}_1},$ where $\boldsymbol{\hat{\prod}}^\dagger =  \boldsymbol{\hat{\prod}}^\ast(\boldsymbol{\hat{\prod}}\boldsymbol{\hat{\prod}}^\ast + \gamma_N\mathbb{I})^{-1}$ and $\boldsymbol{\prod}^\dagger =  \boldsymbol{\prod}^\ast(\boldsymbol{\prod}\boldsymbol{\prod}^\ast + \gamma_N\mathbb{I})^{-1},$ with $(\gamma_N)\subset(0,\infty)$ satisfying $\gamma_N \to 0,$ and with $\boldsymbol{\prod}^\ddagger\coloneqq \boldsymbol{\prod}\boldsymbol{\prod}^\dagger$ and $\Psi'_{[q]} = -B_q\boldsymbol{\prod}$ in \eqref{equation of invertible operators in ARMA(p,q), vector representation}, for $B_q = B_q(M)$ in Assumption \ref{As: projection on finite-dim space for B_q in ARMA(p,q)}, we obtain
\begin{align}
  \hat{B}_q - B_q(M) &= \big(\Psi'_{[q]}  - \hat{\Psi}'_{[q]}\big)\boldsymbol{\hat{\prod}}^\dagger\coprod^{\hat{h}_{\!M}}_{\hat{h}_1} + \,\Psi'_{[q]} \Big(\boldsymbol{\prod}^\dagger\!\coprod^{h_{\!M}}_{h_1} - \boldsymbol{\hat{\prod}}^\dagger\!\coprod^{\hat{h}_{\!M}}_{\hat{h}_1}\Big) + B_q\Big(\boldsymbol{\prod}^\ddagger\!\coprod^{h_{\!M}}_{h_1} - \,\mathbb{I}\Big).\notag
\end{align}
\noindent Further, as in the proof of Lemma \ref{Lemma for ARMA(1,1) operators, finite-dim}, together with \eqref{inequality of difference of estimate to in ARMA(p,q)} and Assumption \ref{As: Dim and eigvals fARMA pq}, it holds
\begin{align*}
  \|\hat{B}_q - B_q(M)\|_\mathcal{S} 
  &\leq (\hat{\zeta}_M + \gamma_N)^{-1/2}\|\hat{\Psi}_{\!L} - \Psi_{\!L}\|_\mathcal{S} + \|\Psi'_{[q]} \|_\mathcal{S}\Big\|\boldsymbol{\prod}^\dagger\!\coprod^{h_{\!M}}_{h_1} - \boldsymbol{\hat{\prod}}^\dagger\!\coprod^{\hat{h}_{\!M}}_{\hat{h}_1}\Big\|_\mathcal{L} + \gamma_N\zeta^{-1}_M\|B_q(M)\|_\mathcal{S}\notag\\
  &\leq \mathrm{O}_{\Prob} \big(M^{1/2}Z^2_M \|\hat{\Psi}_{\!L} - \Psi_{\!L}\|_\mathcal{S} \big) + \mathrm{O}(\gamma_N \zeta^{-1}_{\!M})
\end{align*}
Hence, as $\gamma_N = \mathrm{O}(K^{-\beta}),$ for fixed $M$ holds $\|\hat{B}_q - B_q(M)\|_\mathcal{S}=\mathrm{O}_{\Prob}(K^{-\beta})$ and thus $\|\hat{\beta}_j - \beta_j(M)\|_\mathcal{S} = \mathrm{O}_{\Prob}(K^{-\beta})$ for all $j.$ Further, $\hat{\alpha}_i = \hat{\psi}_i + \hat{B}_{[i]}\hat{\Psi}'''^{\,\top}_{[i]},$ $\alpha_i(M) = \psi_i + B_{[i]}(M)\Psi'''^{\,\top}_{[i]}\!,$ and $\alpha_i = \beta_j=\hat{\beta}_j=0$ for $i>p,$ $j>q,$ give for $L\geq p+q+i-1$ for any $i,p,q\colon$
\begin{align*}
 \|\hat{\alpha}_i-\alpha_i(M)\|_\mathcal{S}\notag &\leq \|\hat{\psi}_i - \psi_i\|_\mathcal{S} + \|\hat{B}_{[i]}(M)\|_\mathcal{L}\|\hat{\Psi}'''_{[i]} - \Psi'''_{[i]}\|_\mathcal{S} + \|\Psi'''_{[i]}\|_\mathcal{L}\|\hat{B}_{[i]} - B_{[i]}(M)\|_\mathcal{S}\notag\\
  &\leq (1 + \|\hat{B}_{[i]}(M)\|_\mathcal{L})\|\hat{\Psi}_{\!L} - \Psi_{\!L}\|_\mathcal{S} + \|\Psi'''_{[i]}\|_\mathcal{L}\|\hat{B}_q - B_q(M)\|_\mathcal{S}.\notag\\
  &\leq \|\hat{B}_q - B_q(M)\|_\mathcal{S}\big(\|\hat{\Psi}_{\!L} - \Psi_{\!L}\|_\mathcal{S} + \|\Psi'''_{[i]}\|_\mathcal{L}\big) + \|\hat{\Psi}_{\!L} - \Psi_{\!L}\|_\mathcal{S},
\end{align*}
and consequently also $\|\hat{\alpha}_i-\alpha_i(M)\|_\mathcal{S}=\mathrm O_{\Prob}(K^{-\beta})$ for all $i.$ This shows our claim.
\end{proof}

\begin{proof}[\textbf{Proof of Theorem \ref{Theo: ARMA(p,q), infinite-dim}}] Follows from the argumentation in the proofs of Theorem \ref{Theorem for ARMA(1,1) operators, infinite-dim} and Lemma \ref{Theorem for ARMA(p,q) operators, finite-dim} and the requirements for our sequences.
\end{proof}

\section*{Acknowledgements}

The authors are grateful to the editor, the  associate editor and the two referees for their constructive comments on an earlier version of this article, which greatly improved the quality of the paper. Parts of the paper were written while Sebastian K\"uhnert was employed at the University of California, Davis. 

\section*{Funding}

Alexander Aue was partially supported by NSF DMS 2515821. Sebastian K\"uhnert was partially supported by TRR 391 \textit{Spatio-temporal Statistics for the Transition of Energy and Transport} (Project number 520388526) funded by the Deutsche Forschungsgemeinschaft (DFG, German Research Foundation). Gregory Rice was partially supported by the Discovery Grant RGPIN 50503-11525 3100 105 from the Natural Science and Engineering Research Council of Canada.

\bibliographystyle{apalike}
\bibliography{References_invertible_lin}

\end{document}